\documentclass[12pt]{amsart}
\usepackage{amsmath,amssymb,bm}
\usepackage{graphicx}
\usepackage{subfigure}
\newtheorem{theorem}{Theorem}[section]
\newtheorem{lemma}[theorem]{Lemma}
\newtheorem{example}[theorem]{Example}
\newtheorem{corollary}[theorem]{Corollary}
\newtheorem{remark}[theorem]{Remark}
\numberwithin{equation}{section}
\numberwithin{table}{section}
\numberwithin{figure}{section}
\allowdisplaybreaks[4]
\begin{document}
\def\cA{\mathcal{A}}
\def\sss{\scriptscriptstyle}
\def\ss{\scriptstyle}
\def\O{{\Omega}}
\def\R{\mathbb{R}}
\def\p{\partial}
\def\cC{\mathcal{C}^{\raise 1pt\hbox{$\ss \PiH$}}}
\def\Ch{\cC_h}
\def\Chk{\cC_{h,k}}
\def\d{\displaystyle}
\def\cT{\mathcal{T}}
\def\KER{K_h^{\raise 1pt\hbox{$\ss \PiH$}}}
\def\AKER{A^{\PiH}_h}
\def\AKERi{A^{\PiH}_{h,i}}
\def\SKER{S^{\PiH}_h}
\def\MKER{M^{\PiH}_{h}}
\def\LMax{\lambda_{\max}(\SKER\AKER)}
\def\LMin{\lambda_{\min}(\SKER\AKER)}
\def\Cond{\kappa(\SKER\AKER)}
\def\Ki{K_{h,i}^{\raise 1pt\hbox{$\ss \PiH$}}}
\def\KERP{\KER\times\KER}
\def\CKer{C_{\mathrm{Ker}\,\Pi}}
\def\MV{{V}_H^{\raise 2pt\hbox{$\ss\rm {ms},h$}}}
\def\MW{{W}_H^{\raise 2pt\hbox{$\ss\rm {ms},h$}}}
\def\MSL{\tV_{H,k}^{\hspace{1pt}\lower 3pt\hbox{$\ss\rm {ms},h$}}}
\def\MSV{V_{H,k}^{\hspace{1pt}\lower 3pt\hbox{$\ss\rm {ms},h$}}}
\def\MSW{W_{H,k}^{\hspace{1pt}\lower 3pt\hbox{$\ss\rm {ms},h$}}}
\def\MSp{{p}_{\sss H}^{\mathrm{ms},h}}
\def\LMSp{{p}_{\sss H,k}^{\mathrm{ms},h}}
\def\MSy{{y}_{\sss H}^{\mathrm{ms},h}}
\def\LMSy{{y}_{\sss H,k}^{\mathrm{ms},h}}
\def\MSzeta{{\zeta}_{\sss H}^{\mathrm{ms},h}}
\def\MSeta{{\eta}_{\sss H}^{\mathrm{ms},h}}
\def\MSpy{(\MSp,\MSy)}
\def\LMSpy{(\LMSp,\LMSy)}
\def\MSpmy{(\MSp,-\MSy)}
\def\MSLu{\tu_{H,k}^{\mathrm{ms},h}}
\def\MSLp{\tp_{H,k}^{\mathrm{ms},h}}
\def\MSLy{\ty_{H,k}^{\mathrm{ms},h}}
\def\MSTwo{\ttV_{H,k}^{\hspace{1pt}\lower 3pt\hbox{$\ss\rm {ms},h$}}}
\def\MSuTwo{\ttu_{H,k}^{\mathrm{ms},h}}
\def\MP{V_{H,\ell }^{\raise 2pt\hbox{$\ss{\rm {ms}},h$}}}
\def\HP{W_{H,\ell }^{\raise 2pt\hbox{$\ss{\rm {ms}},h$}}}
\def\MPu{u_{H,\ell }^{\raise 1.5pt\hbox{$\ss{\rm {ms}},h$}}}
\def\MPp{p_{H,\ell }^{\raise 1.5pt\hbox{$\ss{\rm {ms}},h$}}}
\def\MPy{y_{H,\ell }^{\raise 1.5pt\hbox{$\ss{\rm {ms}},h$}}}
\def\tpsi{\tilde{\psi}}
\def\vT{v_{\sss T}}
\def\bv{\bm{v}}
\def\bw{\bm{w}}
\def\bz{\bm{z}}
\def\bI{\bm{I}}
\def\argmin{\mathop{\rm argmin}}
\def\LT{{L_2(\O)}}
\def\LTP{{\LT\times\LT}}
\def\HO{{H^1(\O)}}
\def\zHO{H^1_0(\O)}
\def\ES{\zHO\times\zHO}
\def\DESh{V_h\times V_h}
\def\DESH{\VH\times \VH}
\def\cB{{\mathcal{B}}}
\def\cV{\mathcal{V}}
\def\cHO{{\mathcal{H}^1}}
\def\cHOP{{\cHO\times\cHO}}
\def\CPF{\mathrm{C_{PF}}}
\def\B#1#2{\cB\big(#1,#2\big)}
\def\py{(p,y)}
\def\pyh{(p_h,y_h)}
\def\pmyh{(p_h,-y_h)}
\def\qz{(q,z)}
\def\tqz{(\tilde{q},\tilde{z})}
\def\qmz{(q,-z)}
\def\tqmz{(\tilde{q},-\tilde{z})}
\def\rs{(r,s)}
\def\rms{(r,-s)}
\def\trs{(\tilde{r},\tilde{s})}
\def\trms{(\tilde{r},-\tilde{s})}
\def\VH{V_{H}}
\def\bB{\mathbb{B}^{\PiH}_h}
\def\bS{\mathbb{S}^{\PiH}_h}
\def\bQ{\mathbf{Q}}
\def\bb{\mathbf{b}}
\def\br{\mathbf{r}}
\def\bp{\mathbf{p}}
\def\bq{\mathbf{q}}
\def\bx{\mathbf{x}}
\def\PiH{\Pi_{\sss H}}
\def\QT{Q_{\sss T}}
\def\tq{\tilde{q}}
\def\tz{\tilde{z}}
\def\bpsi{\bm{\psi}}
\def\bxi{\bm{\xi}}
\def\PFfactor{[1+(\CPF/\alpha)]}
\def\aP{{a\times a}}
%
\title[Multiscale FEMs for an Elliptic Optimal Control Problem]
{Multiscale Finite Element Methods for an Elliptic Optimal Control Problem
 with Rough Coefficients}
 \author{Susanne C. Brenner}
\address{Department of Mathematics and Center for Computation \& Technology,
Louisiana State University, Baton Rouge, LA 70803, USA}\email{brenner@math.lsu.edu}
\author{Jos\'e C. Garay}
\address{Department of Mathematics and Center for Computation \& Technology,
Louisiana State University, Baton Rouge, LA 70803, USA}
\email{jgaray@cct.lsu.edu}
\author{Li-yeng Sung}
\address{Department of Mathematics and Center for Computation \& Technology,
Louisiana State University, Baton Rouge, LA 70803, USA}\email{sung@math.lsu.edu}
\thanks{This work   was supported in part
 by the National Science Foundation under Grant No.
 DMS-19-13035. }
%
%
\begin{abstract}
 We investigate multiscale finite element methods for an elliptic distributed
 optimal control problem with
 rough coefficients.  They are based on the (local)
  orthogonal decomposition methodology of M\aa lqvist and
 Peterseim.
\end{abstract}
\keywords{multiscale, rough coefficients, elliptic optimal control,
localized orthogonal decomposition, domain decomposition}
\subjclass{65N30, 65N15, 65N55, 49N10}
\maketitle
\section{Introduction}\label{sec:Introduction}
 Let $\O$ be a polyhedral domain in $\R^d$ ($d=1, 2,3$) and $y_d\in \LT$.  We consider the following
 elliptic distributed optimal control problem:
\begin{equation}\label{eq:OCP}
 \text{Find}\quad (\bar y,\bar u)= \argmin_{(y,u)\in K}\frac12
  \big[\|y-y_d\|_\LT^2+\gamma\|u\|_\LT^2\big],
\end{equation}
 where $(y,u)\in \zHO\times\LT$ belongs to $K$ if and only if
\begin{equation}\label{eq:PDEConstraint}
     a(y,z)=\int_\O uz\,dx\qquad\forall\,z\in \zHO,
\end{equation}
 and the bilinear form $a(\cdot,\cdot)$ is given by
\begin{equation}\label{eq:aDef}
  a(y,z)=\int_\O (\cA\nabla y)\cdot\nabla z\,dx.
\end{equation}
 Here $\bar y$ is the optimal state, $\bar u$ is the optimal control and $y_d$ is the desired state.
\begin{remark}\label{rem:Notation} \rm
  We will follow the standard notation for differential operators, function spaces and norms
  that can be found for example in \cite{Ciarlet:1978:FEM,BScott:2008:FEM,ADAMS:2003:Sobolev}.
\end{remark}
\par
 We assume only that the components of the symmetric  diffusion matrix $\cA$
  belong to $L_\infty(\O)$
 and the eigenvalues of $\cA$ are bounded below (resp., above) by the positive
  constant $\alpha$ (resp., $\beta$),
 which covers many multiscale optimal control problems.
\begin{example}\label{example:Oscillatory}\rm
  This example is from \cite{HW:1999:MS}, where $\O$ is the unit square $(0,1)\times (0,1)$,
\begin{equation*}
 \cA(x)= \begin{bmatrix}
    c(x) &0 \\ 0 &c(x)
  \end{bmatrix},
\end{equation*}
 and
\begin{equation*}
   c(x)=\frac{\d 2+1.8\sin\Big(\frac{2\pi x_1}{\epsilon}\Big)}
   {\d 2+1.8\cos\Big(\frac{2\pi x_2}{\epsilon}\Big)}
     +\frac{\d 2+\sin\Big(\frac{2\pi x_2}{\epsilon}\Big)}
     {\d 2+1.8\sin\Big(\frac{2\pi x_1}{\epsilon}\Big)}
\end{equation*}
 is highly oscillatory for small $\epsilon$.  Note that
 $$\alpha=\min_{0\leq x\leq 1}c(x)\approx 1.248\quad\text{and}\quad
    \beta=\max_{0\leq x\leq 1}c(x)\approx 19.526$$
  for any $\epsilon\leq 1$.
\end{example}
\begin{example}\label{example:Heterogeneous}\rm
  This example is from \cite{BGS:2021:LOD}, where $\O$ is the unit square $(0,1)\times(0,1)$,
 $$\cA=\begin{bmatrix}
   \cA_{11}&0\\ 0 &\cA_{22}
 \end{bmatrix},$$
 and the components $\cA_{11}$ and $\cA_{22}$ are randomly generated
 piecewise constant functions with respect
 to a uniform partition of $\O$ into $40\times 40$ small squares
 (cf. Figure~\ref{fig:Heterogeneous}).
 The values of $\cA_{11}$ and $\cA_{22}$ are between $\alpha=1$ and $\beta=1350$.
\begin{figure}[hhh]
\begin{center}
  \includegraphics[width=.3\linewidth]{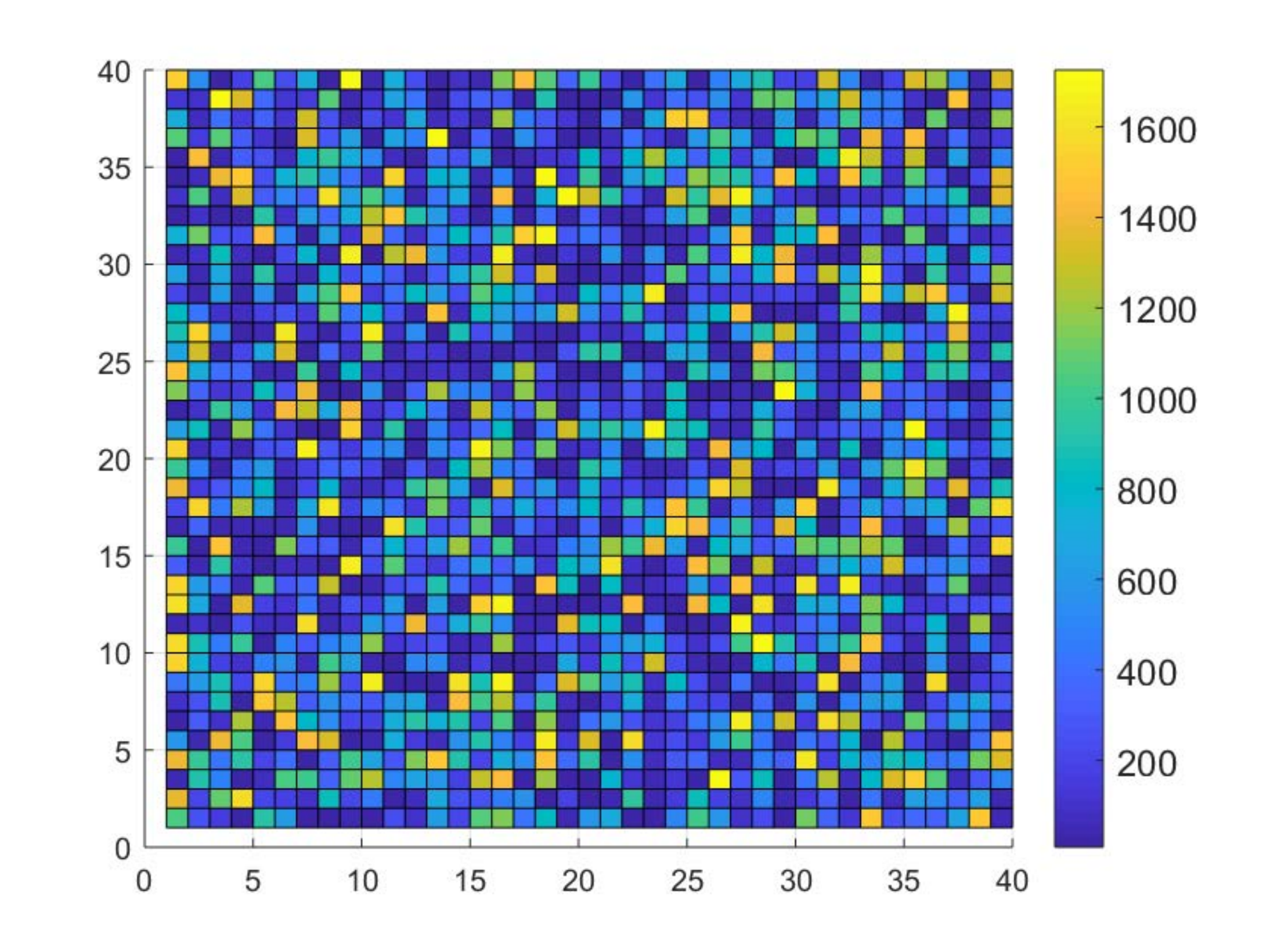} \hspace{20pt}
  \includegraphics[width=.3\linewidth]{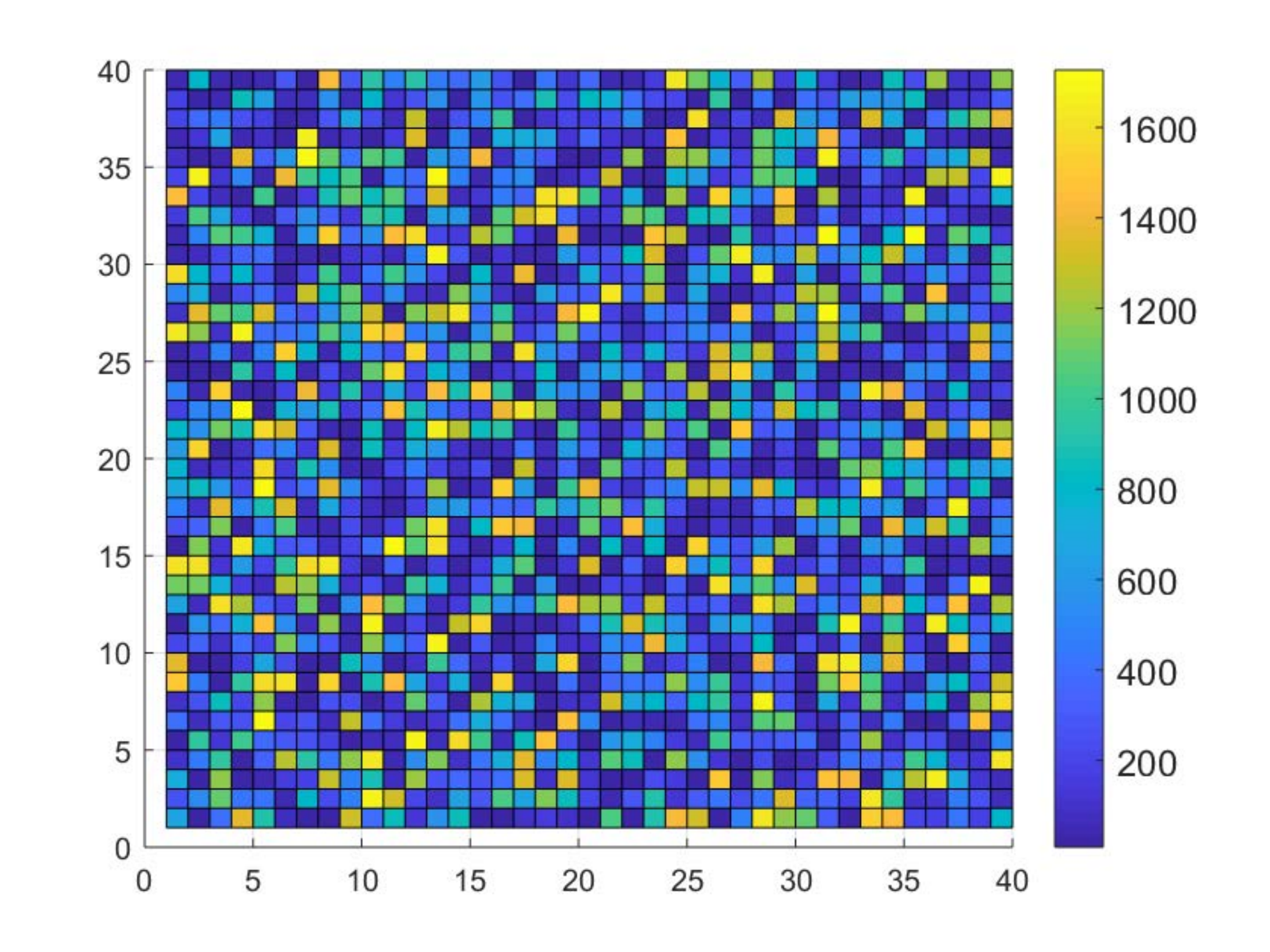}
\end{center}
\caption{$\cA_{11}$ $($left$)$ and $\cA_{22}$ $($right$)$}
\label{fig:Heterogeneous}
\end{figure}
\end{example}
\par
 Due to the roughness of the coefficients in \eqref{eq:aDef}, a
 standard finite element method can only accurately
 capture the
 solution of \eqref{eq:OCP} on a very fine mesh (cf. \cite{BO:2000:Bad}),
 which can be too expensive,
 especially when
 the problem has to be solved repeatedly for different $y_d$.
 Our goal is to construct generalized finite element spaces that can
 produce approximate solutions of
 \eqref{eq:OCP} with $O(H)$ (resp., $O(H^2)$) error in the energy
 (resp., $L_2$) norm, where $H$ is the
 mesh size and the dimensions of the generalized finite element spaces are
  $O(H^{-d})$.  In other words
 the performance of these generalized finite element methods is similar
  to standard finite element methods
 for elliptic problems with smooth coefficients on smooth or convex domains.
\par
  Our constructions are based on the Localized Orthogonal Decomposition
   (LOD) approach in
  \cite{HP:2013:LOD,MP:2014:LOD} and the ideas in
  \cite{KPY:2018:LOD,BGS:2021:LOD} (cf. also
  \cite[Section~4.3]{MP:2021:LOD}).  The basis functions of the
  generalized finite element spaces are obtained by a correction process
  that can be carried out offline.
  The online computation only involves solving a  linear system of moderate size.
  Therefore these generalized finite
  element methods can also be viewed as reduced order methods that are particularly
  suitable for repeat solves.
\par
  There are many numerical methods for elliptic problems with rough coefficients besides the
  LOD methods.  They include  the variational multiscale method
  (cf. \cite{Hughes:1995:Multiscale,HFMQ:1998:VMS,HS:2007:VM} and the references therein),
  the multiscale finite element method
 (cf. \cite{HW:1999:MS,HWC:1999:Multiscale,EH:2009:MSFEM,HOS:2014:AdaptiveMSFEM}
 and the references therein), the heterogeneous multiscale method
 (cf. \cite{EE:2003:HMSFEM,EMZ:2005:HMS,Abdulle:2005:HMS,AEEV:2012:HMS}
  and the references therein), and the method of approximate component synthesis
 (cf. \cite{HL:2010:Special,HK:2014:Special} and the references therein).
  We refer the readers to \cite{OS:2019:Book,MP:2021:LOD} for the discussion of other methods.
\par
 On the other hand, as far as we know, there is only one paper \cite{GYWLY:2018:HMMOCP} that solved the
 optimal control problem \eqref{eq:OCP}--\eqref{eq:aDef} (with additional control constraints)
 by the heterogeneous multiscale
 method, where scale separation and periodic structure are assumed.
 In the context of parabolic optimal control problems with rough coefficients, reduced
 order finite element methods in the same spirit of the current paper are treated in
 \cite{MPS:2018:LODParabolic,ZCL:2020:Parabolic}.  In particular, the methodology in
 \cite{MPS:2018:LODParabolic} is also based on the LOD approach.  The distinctive
 feature of our work in this paper is that the construction of the localized multiscale finite
  element space and its analysis are based entirely  on classical techniques from
  domain decomposition and numerical linear algebra.
\par
 The rest of the paper is organized as follows.  We recall relevant results for the optimal control
 problem in Section~\ref{sec:Continuous}.  A multiscale finite element method based on
 orthogonal decomposition is treated in Section~\ref{sec:Ideal}.
 We introduce a localized multiscale finite element space in Section~\ref{sec:LMFES} and analyze
 the
 corresponding Galerkin method in Section~\ref{sec:LMFEM}, where the error estimates in
 Section~\ref{sec:Ideal} play a  useful role.
 Numerical results are presented in Section~\ref{sec:Numerics} and we end with some concluding
 remarks in Section~\ref{sec:Conclusion}.
\par
 We will use $\langle\cdot,\cdot\rangle$ to denote the canonical bilinear form on a
   finite dimensional vector space $V$ and its dual
 space $V'$.
 A linear  operator $L:V\longrightarrow V'$ is symmetric if
\begin{equation*}
 \langle Lv_1,v_2\rangle=\langle Lv_2,v_1\rangle \qquad \forall\,v_1,v_2\in V,
\end{equation*}
 and it is symmetric positive definite (SPD) if additionally
\begin{equation*}
 \langle Lv,v\rangle>0 \qquad\forall\,v\in V\setminus\{0\}.
\end{equation*}
 Given two finite dimensional vector spaces $V$ and $W$ and a linear transform
   $T:V\longrightarrow W$,
 the transpose $T^t:W'\longrightarrow V'$ is defined by
\begin{equation*}
 \langle  T^t\mu,v\rangle=\langle \mu,Tv\rangle \qquad\forall\,\mu\in W',\,v\in V.
\end{equation*}
\par
 We also assume that all the unspecified positive constants in the paper are greater than or equal to 1.
\section{The Continuous Problem}\label{sec:Continuous}
 By a standard result \cite[Section~2.2]{Lions:1971:OC},
 the convex minimization problem defined by \eqref{eq:OCP}--\eqref{eq:aDef}
  has a unique solution
 determined by  the following first order optimality conditions:
 \begin{alignat*}{3}
  a(\bar y,z)&=\int_\O \bar uz\,dx&\qquad&\forall\,z\in \zHO,\\
  a(\bar p,q)&=\int_\O (y_d-\bar y)q\,dx &\qquad& \forall\,q\in \zHO,\\
  \bar p&=\gamma \bar u,
 \end{alignat*}
  where the adjoint state $\bar p$ belongs to $\zHO$.
\par
 After eliminating $\bar u$,  we have  the following system for $(\bar p,\bar y)$:
\begin{alignat}{3}
    a(\bar p,q)+\int_\O \bar y q\,dx&=\int_\O y_d q\,dx&\qquad&\forall\,q\in\zHO,\label{eq:C1}\\
    \int_\O \bar p z\,dx-\gamma a(\bar y,z)&=0&\qquad&\forall\,z\in\zHO.\label{eq:C2}
\end{alignat}
\begin{remark}\label{rem:ContinuousProblem}\rm
 Note that \eqref{eq:C1}--\eqref{eq:C2} is equivalent to
\begin{alignat*}{3}
    \tilde a(\bar p,q)+\int_\O \tilde y q\,dx&=\int_\O \tilde y_d q\,dx&\qquad&\forall\,q\in\zHO,\\
    \int_\O \bar p z\,dx-\tilde\gamma \tilde a(\tilde y,z)&=0&\qquad&\forall\,z\in\zHO,
\end{alignat*}
 where $\tilde a(\cdot,\cdot)=\tau a(\cdot,\cdot)$,
  $\tilde y=\tau\bar y$,
 $\tilde y_d=\tau y_d$, $\tilde\gamma=(\gamma/\tau^2)$ and
  $\tau$ is any positive number.
 Therefore we can assume that the lower bound $\alpha$ for the eigenvalues of
 $\cA$ (cf. \eqref{eq:aDef}) in the definition of the  bilinear form $a(\cdot,\cdot)$
 in \eqref{eq:C1}--\eqref{eq:C2}
 is roughly $1$, as in Example~\ref{example:Oscillatory} and Example~\ref{example:Heterogeneous}.
\end{remark}
\par
  Since the dependence on $\gamma$  is not our main concern here, we will
 take $\gamma$ to be 1 in \eqref{eq:C2}.
 We will also simplify the notation by dropping the bars over $p$ and $y$
 and
 consider the problem of finding $(p,y)\in\zHO\times\zHO$ such that
\begin{alignat}{3}
    a(p,q)+\int_\O y q\,dx&=\int_\O y_d q\,dx&\qquad&\forall\,q\in\zHO,
    \label{eq:KKT1}\\
    \int_\O p z\,dx-a(y,z)&=0&\qquad&\forall\,z\in\zHO.
    \label{eq:KKT2}
\end{alignat}
\par
 We can write \eqref{eq:KKT1}--\eqref{eq:KKT2} concisely as
\begin{equation}\label{eq:Concise}
   \B{\py}{\qz}=\int_\O y_d q\,dx\qquad\forall\,\qz\in \ES,
\end{equation}
 where
\begin{equation}\label{eq:cBDef}
  \B{\py}{\qz}=a(p,q)+\int_\O yq\,dx+\int_\O pz\,dx-a(y,z).
\end{equation}
\par
  We will use $\|\cdot\|_a$ to denote the energy norm
 $\sqrt{a(\cdot,\cdot)}$.
 Note that
\begin{equation}\label{eq:Comparison}
  \sqrt\alpha\,|v|_\HO\leq \|v\|_a \leq\sqrt\beta\,|v|_\HO \qquad
   \forall\,v\in \HO
\end{equation}
 by our assumption on $\cA$, and we have a Poincar\'e-Friedrichs inequality
 \cite{ADAMS:2003:Sobolev}
\begin{equation}\label{eq:PF}
  \|v\|_\LT^2\leq \CPF |v|_\HO^2\qquad\forall\,v\in\zHO.
\end{equation}
\par
 The following are the salient features of $\B{\cdot}{\cdot}$ that
 follow immediately from \eqref{eq:cBDef}--\eqref{eq:PF} and the
 Cauchy-Schwarz inequality:
\begin{equation}\label{eq:mCoercive}
   \B{\qz}{\qmz}= \|\qz\|_{a\times a}^2\quad\forall\,\qz\in\ES,
\end{equation}
 and
 \begin{equation} \label{eq:Bounded}
  \B{\qz}{\rs}\leq \PFfactor\|\qz\|_\aP\|\rs\|_\aP
\end{equation}
 for all $\qz, \rs\in\ES$,
 where the norm $\|\cdot\|_{a\times a}$ is defined by
\begin{align}
\|\qz\|_{a\times a}^2&=\|q\|_a^2+\|z\|_a^2.\label{eq:ProdEnergyNoem}
\end{align}
\par
 From here on we will also use bold-faced letters to denote
 members of the product space $\ES$ in order to improve the
 readability of the formulas.
\begin{lemma}\label{lem:BStability}
 Let $V$ be a subspace of $H^1_0(\O)$.
  We have
\begin{equation}\label{eq:InfSup}
   \inf_{\strut\bv\in V\times V}\sup_{\substack{\vspace{1pt}\\ \bw\in V\times V}}
   \frac{\B{\bv}{\bw}}{\|\bv\|_\aP\|\bw\|_\aP}\geq 1.
\end{equation}
\end{lemma}
\begin{proof}  Let $\bv=\qz\in V\times V$ be arbitrary.
 According to \eqref{eq:mCoercive},  we have
\begin{equation*}
  \|\qz\|_\aP^2=\B{\qz}{\qmz}
\end{equation*}
 and consequently
\begin{equation}\label{eq:BS1}
  \|\bv\|_\aP= \frac{\B{\qz}{\qmz}}{\|\qz\|_\aP}
        =\frac{\B{\qz}{\qmz}}{\|\qmz\|_\aP}
        \leq \sup_{\bw\in V\times V}\frac{\B{\bv}{\bw}}{\|\bw\|_\aP}.
\end{equation}
\end{proof}
\begin{remark}\label{rem:Isomorphism}\rm
 Let $V$ be a closed subspace of $H^1_0(\O)$ and $\bm{V}=V\times V$.  It follows from
 \eqref{eq:Bounded} that we can define a linear transformation $T:\bm{V}\longrightarrow \bm{V}'$
 by
\begin{equation*}
  \langle T\bz,\bw\rangle=\B{\bz}{\bw} \qquad\forall\,\bz,\bw\in \bm{V}.
\end{equation*}
 Since the bilinear form $\B{\cdot}{\cdot}$ is symmetric, the
  inf-sup condition \eqref{eq:InfSup}
 implies that $T$ is an isomorphism
 and the operator norms of $T$ and $T^{-1}$ (with respect to $\|\cdot\|_\aP$)
 are bounded by $1$
 (cf. \cite{Babushka:1973:LM,Brezzi:1974:SPP}).
\end{remark}
\begin{remark}\label{rem:FEMs}\rm
 In view of  Remark~\ref{rem:Isomorphism}, one can solve
 \eqref{eq:Concise} by a standard finite element method.
 Let $V_h\subset \zHO$ (resp., $\VH\subset \zHO$) be the $P_1$ or $Q_1$
 finite element space
  associated with the triangulation $\cT_h$ (resp., $\cT_H$) of $\O$,
  where $\cT_h$ is a refinement of $\cT_H$ and hence $V_H$ is a subspace of $V_h$.
\par
 We assume that $h\ll 1$ so that $\pyh\in \DESh$ determined by
\begin{equation}\label{eq:FineDiscreteProblem}
  \B{(p_h,y_h)}{\qz}=\int_\O y_d q\,dx\qquad\forall \,\qz\in \DESh
\end{equation}
 provides a good approximation of the solution $\py$ of \eqref{eq:Concise}, but the
 dimension of $V_h$ is so large that the computational cost is prohibitive,
 especially if we
 have to solve \eqref{eq:FineDiscreteProblem} repeatedly for different $y_d$.
\par
 On the other hand, for $H\gg h$, the solution
 $(p_{\sss H},y_{\sss H})\in \DESH$ defined by
\begin{equation}
  \B{(p_{\sss H},y_{\sss H})}{\qz}=\int_\O y_d q\,dx\qquad\forall \,\qz\in \DESH
\end{equation}
 is computationally feasible but
 not sufficiently accurate.  Therefore we need generalized finite element spaces
 to bridge the two scales.
\end{remark}
\par
\begin{figure}[htbp]
\begin{center}
  \includegraphics[width=.3\linewidth]{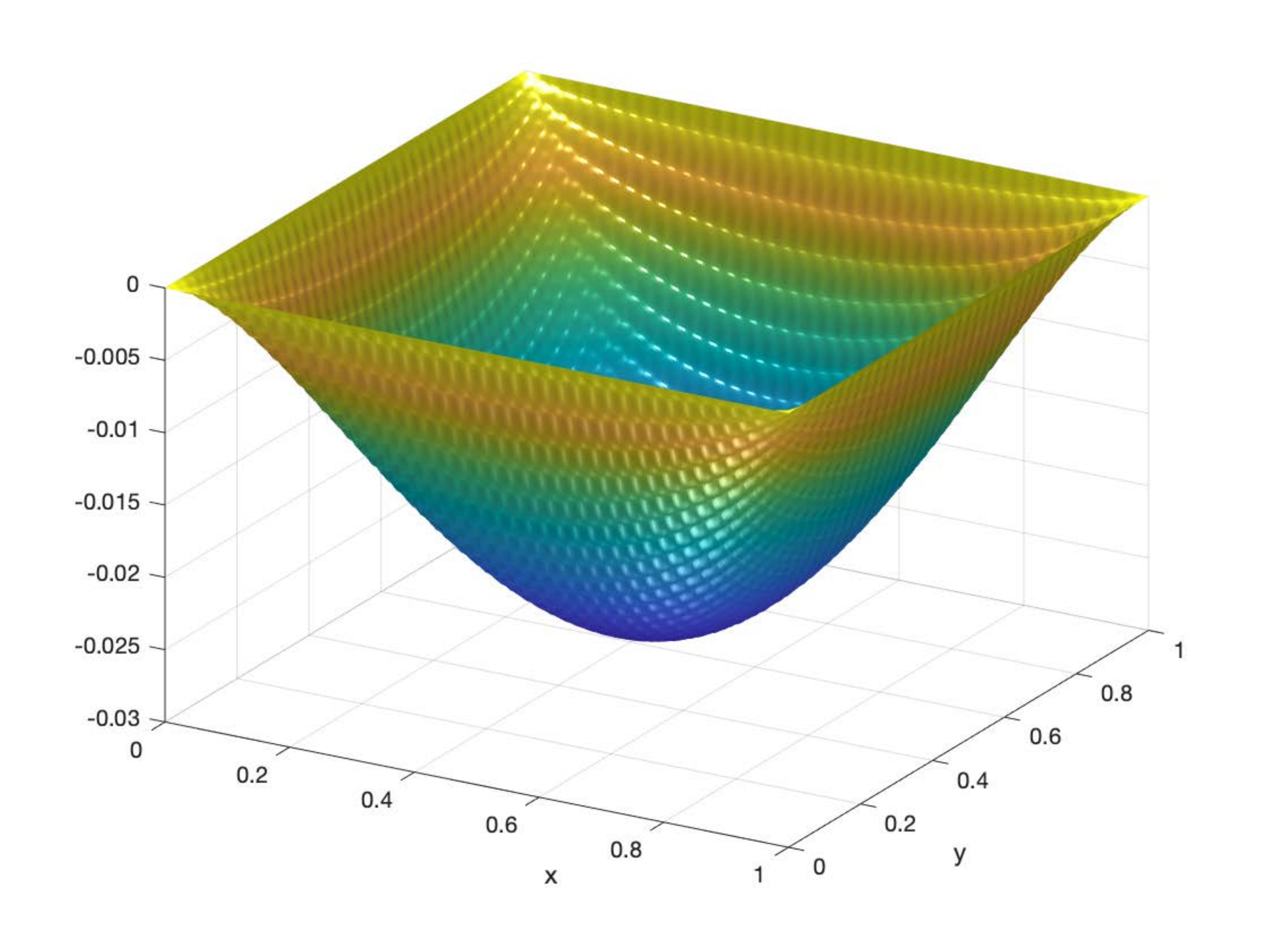} \hspace{20pt}
  \includegraphics[width=.3\linewidth]{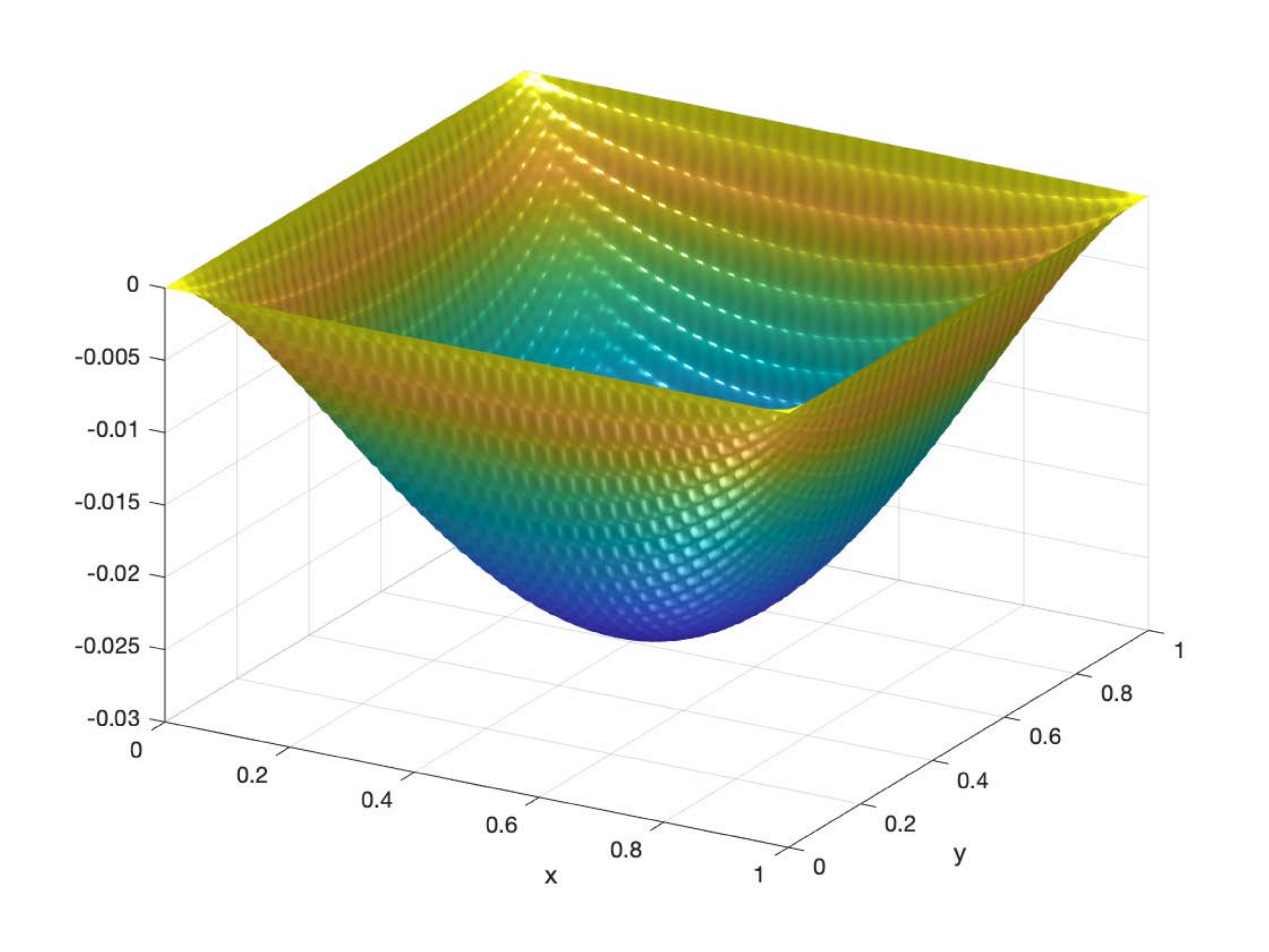}
\par
   \includegraphics[width=.3\linewidth]{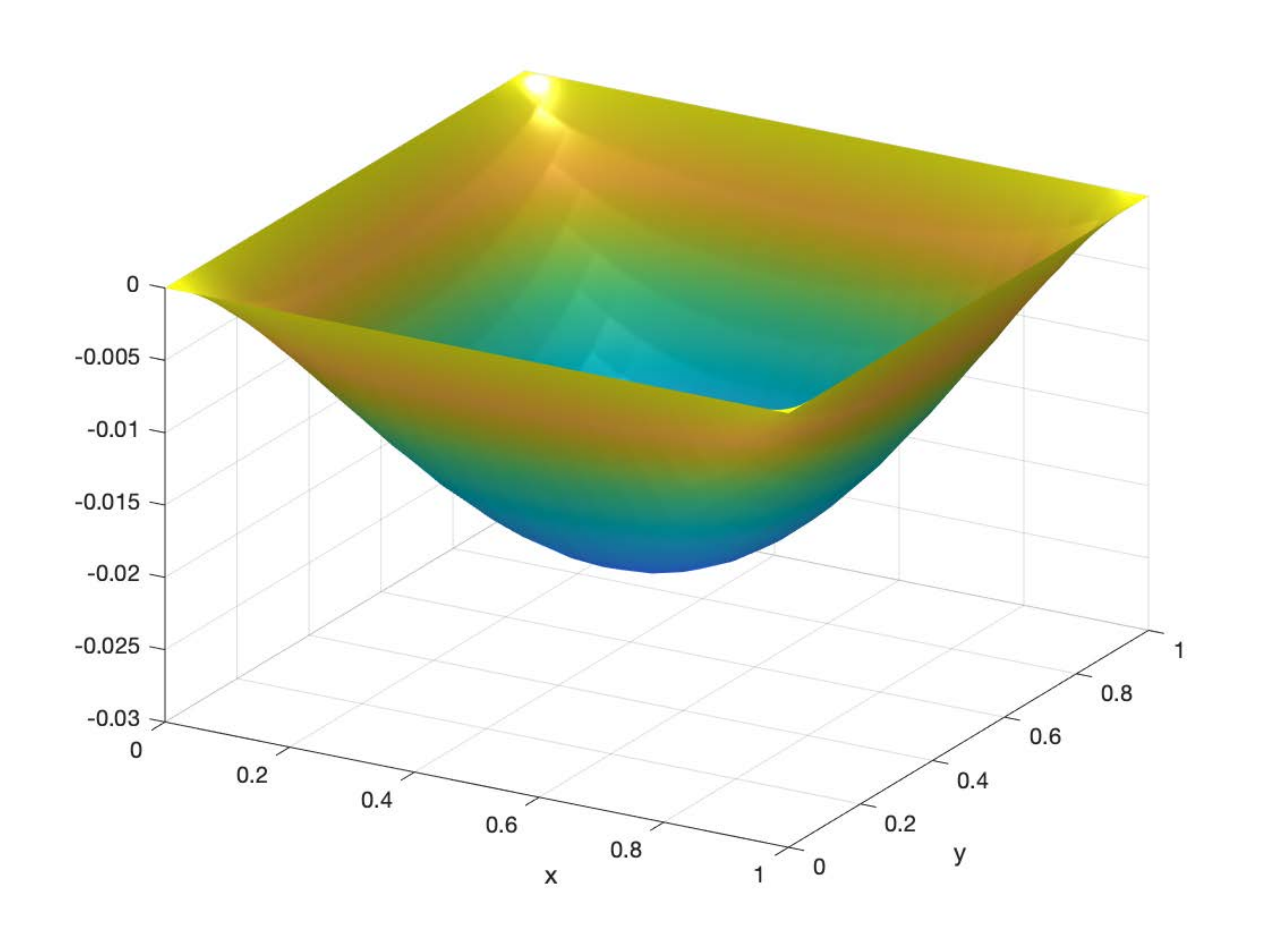}
\end{center}
\caption{Finite element solutions of the optimal state in Example~\ref{example:Oscillatory}
with $\epsilon=0.025$:
 solution on a fine mesh with $h=1/320$ (top left), solution on a coarse mesh with $H=1/20$
 (bottom), LOD solution with $H=1/20$ and $h=1/320$ (top right)}
\label{fig:OscillatoryFEMs}
\end{figure}
\begin{figure}[hhhh]
\begin{center}
  \includegraphics[width=.3\linewidth]{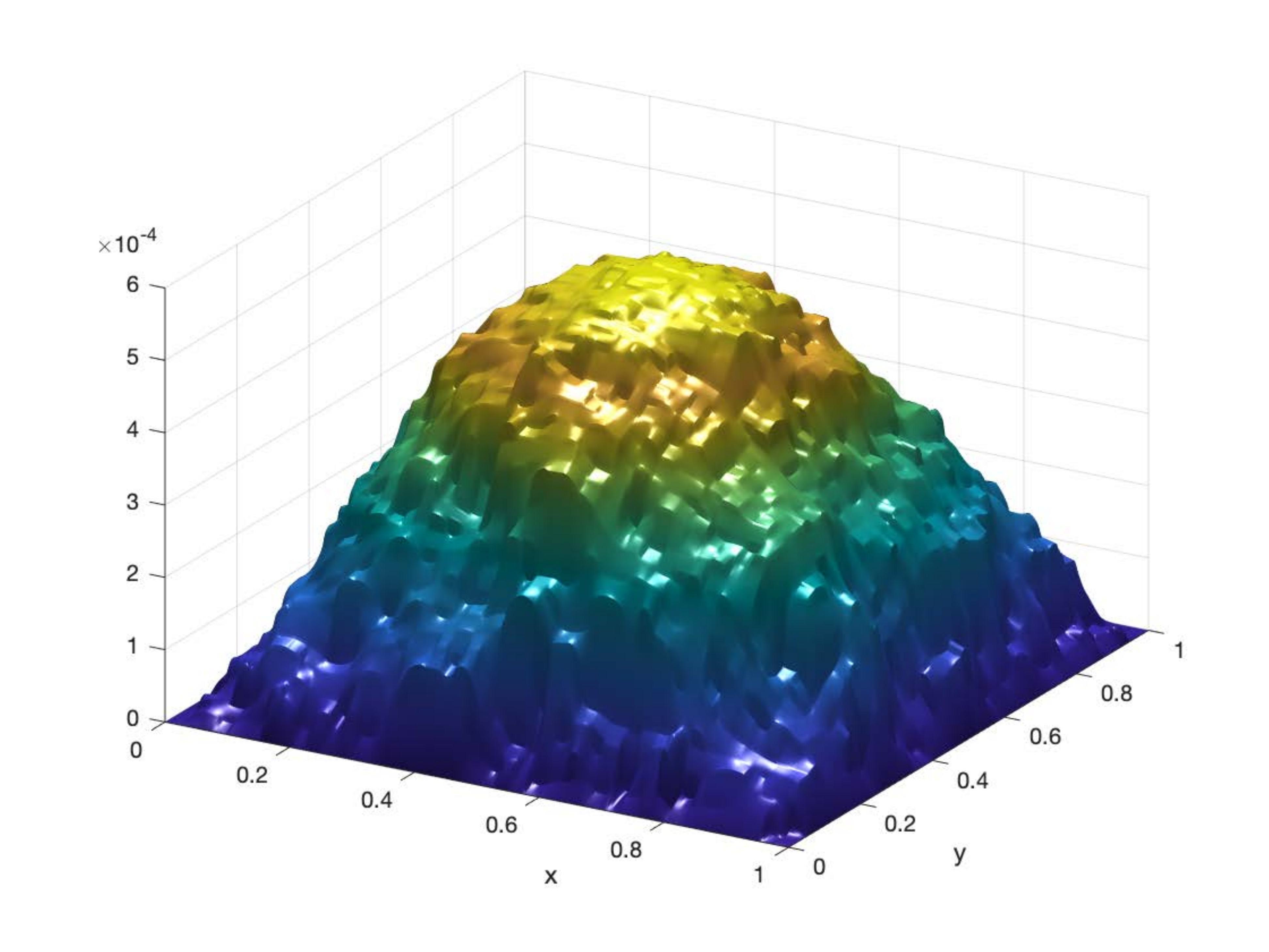} \hspace{20pt}
  \includegraphics[width=.3\linewidth]{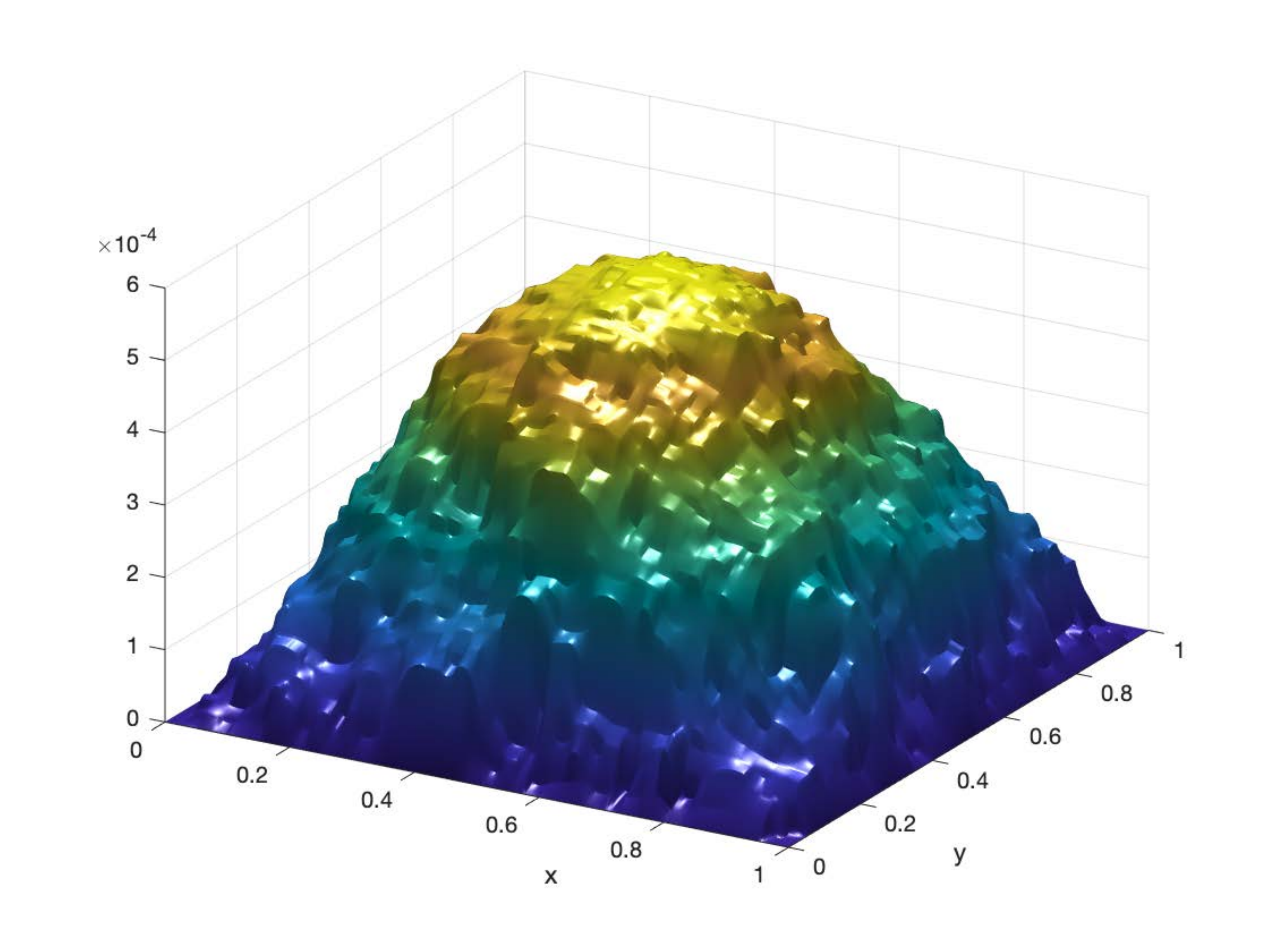}
\par
   \includegraphics[width=.3\linewidth]{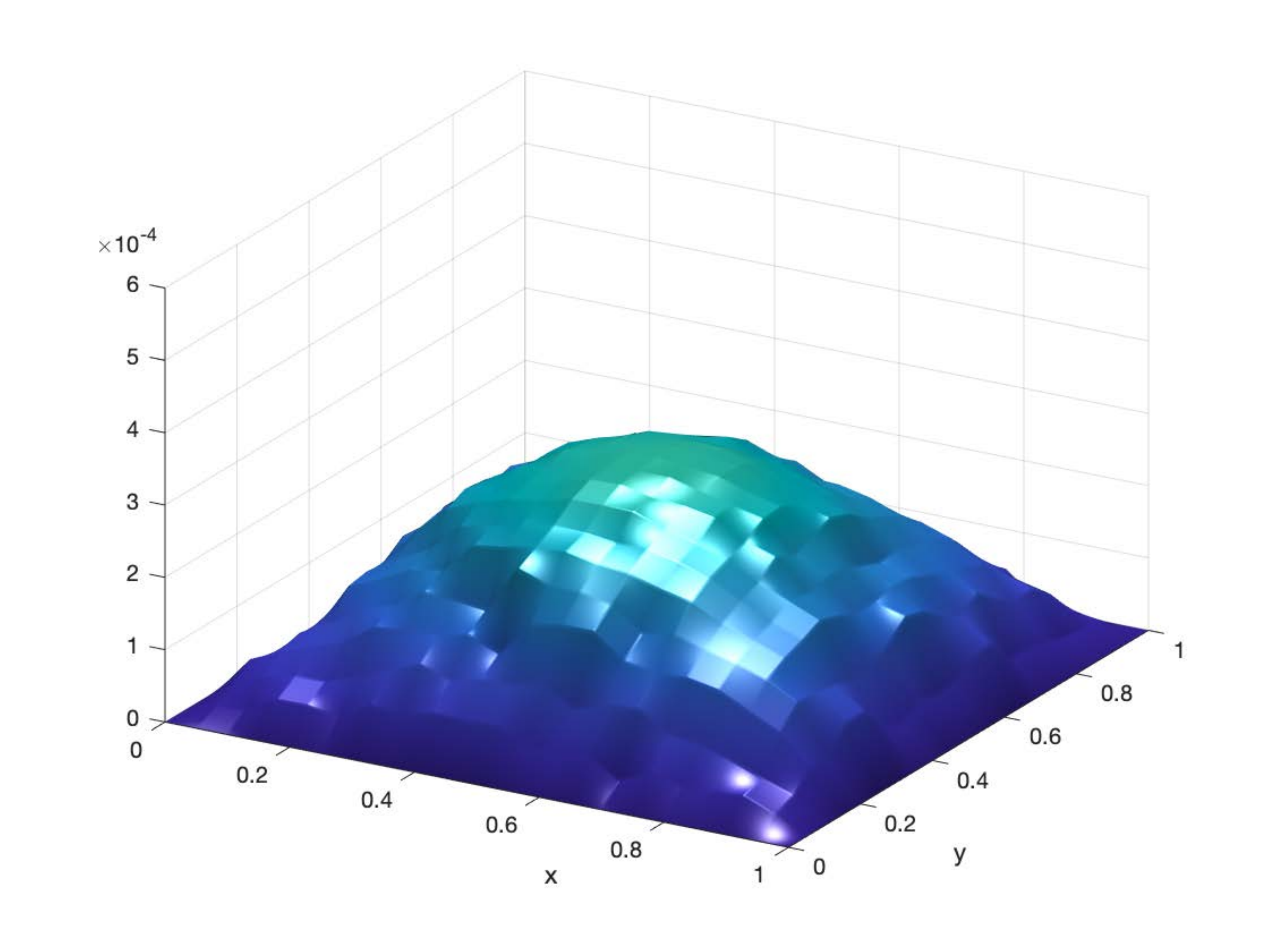}
\end{center}
\caption{Finite element solutions of the optimal state in Example~\ref{example:Heterogeneous}:
 solution on a fine mesh with $h=1/320$ (top left), solution on a coarse mesh with $H=1/20$
 (bottom), LOD solution with $H=1/20$ and $h=1/320$ (top right)}
\label{fig:HeterogenousFEMs}
\end{figure}

  Finite element solutions for the optimal states in Example~\ref{example:Oscillatory} and
  Example~\ref{example:Heterogeneous} are displayed in Figure~\ref{fig:OscillatoryFEMs} and
  Figure~\ref{fig:HeterogenousFEMs}.  It can be observed for both examples
   that the LOD solutions from Section~\ref{sec:LMFEM} capture the
  fine scale solutions while the coarse scale solutions fail to do so.
\begin{remark}\label{rem:WPBound}\rm
  It follows from  \eqref{eq:Comparison}, \eqref{eq:PF},  \eqref{eq:InfSup}
  and \eqref{eq:FineDiscreteProblem} that
\begin{align*}
  \|\pyh\|_\aP\leq \|y_d\|_\LT\sup_{\qz\in\DESh} \frac{\|q\|_\LT}{\|\qz\|_\aP}
  \leq \sqrt{\CPF/\alpha}\,\|y_d\|_\LT.
\end{align*}
\end{remark}
\section{The Ideal Multiscale Finite Element Method}\label{sec:Ideal}
\par
 In this Section we construct and analyze the ideal multiscale  finite element method
  following the ideas in
 \cite{MP:2014:LOD,KPY:2018:LOD}, which begins with the construction of
 a projection operator.
 We will denote by $n$ (resp., $m$) the dimension of the finite element
  space $V_h$ (resp., $V_H$)
 in Remark~\ref{rem:FEMs}.
\subsection{The Projection Operator $\PiH$}
 The operator $\PiH:H^1_0(\O)\longrightarrow V_H$ is defined by taking
  the nodal average of the local $L_2$
    orthogonal projections of $\zeta\in H^1_0(\O)$ into $P_1$ or $Q_1$ polynomials.
    More precisely,  we define
    $\PiH\zeta$ by
\begin{equation}\label{eq:PiHDef}
  (\PiH \zeta) (p)=\frac{1}{|\cT_p|}\sum_{T\in\cT_p} (\QT \zeta_{\sss T})(p)
  \qquad\forall\,p\in\cV_H,
\end{equation}
 where $\cV_H$ is the set of all the (interior) vertices of $\cT_H$, $\cT_p$
 is the set of the elements in
 $\cT_H$ that share $p$ as a common vertex, $|\cT_p|$ is the number of
 elements in $\cT_p$,
 $\zeta_{\sss T}$ is the restriction of $\zeta$ to $T$, and
 $\QT$ is the orthogonal projection from $L_2(T)$ onto $P_1(T)$ or $Q_1(T)$.
\par
 We have an obvious relation
\begin{equation}\label{eq:Projection}
   \PiH v=v \qquad\forall\,v\in V_H
\end{equation}
 and also an interpolation error estimate \cite[Appendix~A]{BGS:2021:LOD}
\begin{equation}\label{eq:PiEstimates}
   H^{-1}\|v-\PiH v\|_\LT
   +|\PiH v|_\HO\leq C_\dag |v|_\HO \qquad \forall\,v\in H^1_0(\O),
\end{equation}
 where
 the positive constant $C_\dag$ depends only on the shape regularity of $\cT_H$.
\begin{remark}\label{rem:AlternativePiHEst}\rm
  We can use \eqref{eq:Comparison} to translate the estimate for $|\PiH v|_\HO$ into
\begin{equation*}
  \|\PiH v\|_a\leq C_\dag\sqrt{\beta/\alpha}\|v\|_a \qquad\forall\,v\in\zHO.
\end{equation*}
\end{remark}
%
\par
\begin{figure}[htbp]
\begin{center}
  \includegraphics[width=1.8in]{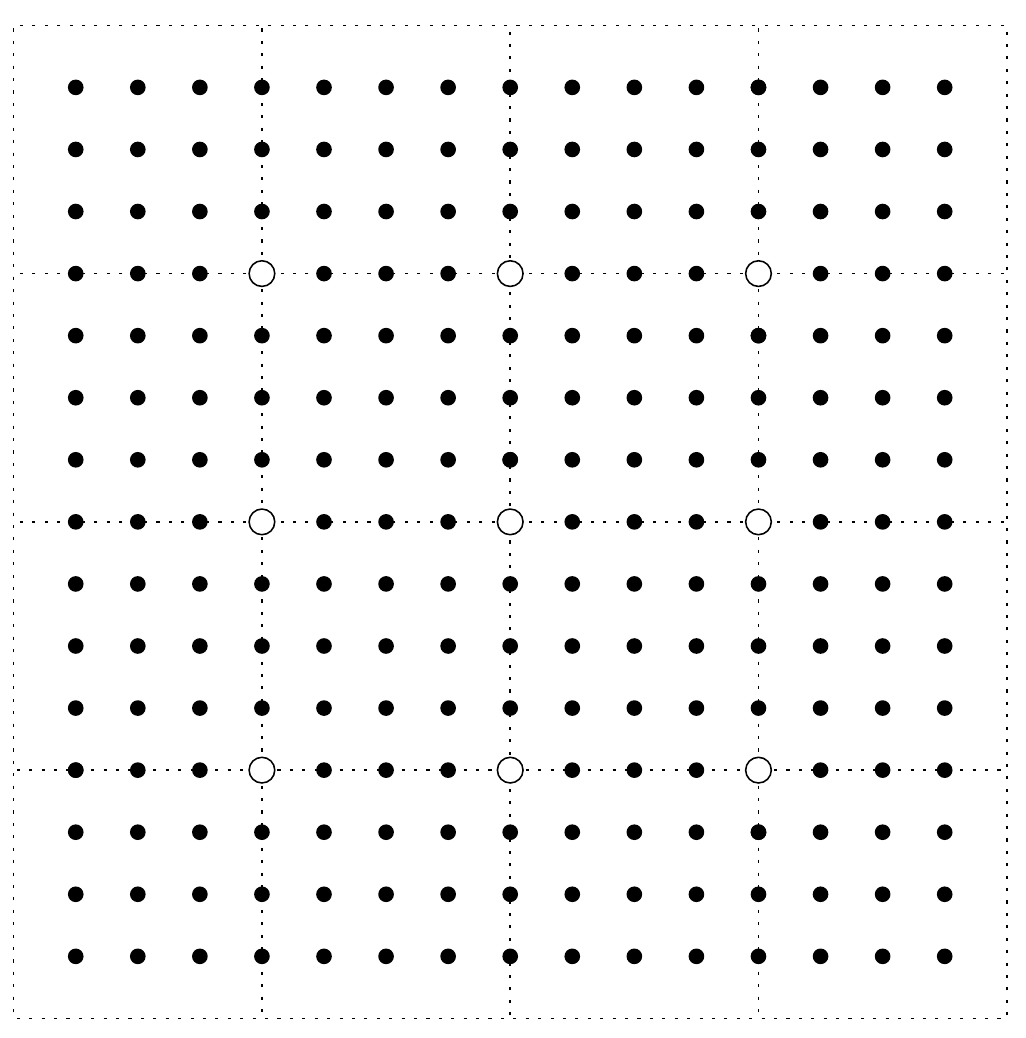}
\end{center}
\caption{The nodes for $V_H$
  are represented by the circles
 and the nodes for $\varphi_1,\ldots,\varphi_\ell\in V_h$ are
 represented by the solid dots.}
\label{fig:Basis}
\end{figure}

 We will denote the kernel of the restriction of $\Pi_H$ to $V_h$ by $\KER$, i.e.,
\begin{equation}\label{eq:KERDef}
 \KER=\{v\in V_h:\,\PiH v=0\}.
\end{equation}
 It follows from \eqref{eq:Projection} and \eqref{eq:KERDef} that
\begin{equation}\label{eq:dimKER}
  \mathrm{dim}\KER=\mathrm{dim}V_h-\mathrm{dim}V_H=n-m.
\end{equation}
 A basis for $\KER$ is given in the following lemma.
\begin{lemma}\label{lem:KERBasis}
 Let $\ell=n-m$  and $\varphi_1,\ldots,\varphi_\ell$ be the nodal basis
 functions in $V_h$ that vanish at the nodes of $V_H$
  $($cf. Figure~\ref{fig:Basis} for a two dimensional example
 with the $Q_1$ finite element$)$.  Then
 $(I-\PiH)\varphi_1,\ldots,(I-\PiH)\varphi_\ell$ form a basis
 of $\KER$, where $I$ is the identity operator on $V_h$.
\end{lemma}
\begin{proof} It follows from \eqref{eq:Projection} that
$(I-\PiH)\varphi_i\in \KER$ for $1\leq i\leq \ell$. In view of
\eqref{eq:dimKER}, it only remains to show that the functions
$(I-\PiH)\varphi_1,\ldots,(I-\PiH)\varphi_\ell$ are linearly independent.
\par
 Suppose $\sum_{i=1}^\ell c_i (I-\PiH)\varphi_i=0$.
 Then the function $\sum_{i=1}^\ell c_i\varphi_i=\sum_{i=1}^\ell c_i\PiH\varphi$
 belongs to $V_H$ and at the same time vanishes at the nodes of $V_H$.  It follows that
  $\sum_{i=1}^\ell c_i\varphi_i=0$ and hence $c_i=0$ for $1\leq i\leq\ell$
  because the functions
  $\varphi_1,\ldots,\varphi_\ell$ are linearly independent.
\end{proof}
%
\subsection{The Projection Operator $\Ch$}\label{subsec:Ch}
 According to Remark~\ref{rem:Isomorphism}, we can
 define a linear transformation $$\Ch:\DESh\longrightarrow\KERP$$ by
\begin{equation}\label{eq:ChDef}
  \B{\Ch\bv}{\bw}=\B{\bv}{\bw}\qquad\forall\,\bv\in \DESh,\;\bw\in\KERP.
\end{equation}
 The elementary algebraic properties of $\Ch$ that follow directly from
 \eqref{eq:ChDef}
 are collected in the following lemma.
\begin{lemma}\label{lem:Ch}
   We have
\begin{alignat}{3}
\B{\Ch\bv}{\bw}&=\B{\bv}{\Ch\bw}&\qquad&\forall\,\bv,\bw\in\DESh,\label{eq:Ch5}\\
(\PiH\times\PiH)\Ch\bv&=0&\qquad&\forall\,\bv\in\DESh,\label{eq:Ch1}\\
\Ch \bv &=\bv &\qquad&\forall\,\bv\in \KERP. \label{eq:Ch2}
\end{alignat}
\end{lemma}
\begin{remark}\label{rem:Projections}\rm
 It follows from \eqref{eq:Ch2} that $\Ch$ is a projection from
 $\DESh$ onto $\KERP$, and that
\begin{equation}\label{eq:IMCh}
 (\bm{I}-\Ch)(\bm{I}-\PiH\times\PiH)\bv=0 \qquad\forall\,\bv\in\DESh,
\end{equation}
 where $\bm{I}$ is the identity operator on $\DESh$.
\end{remark}
\par
\begin{lemma}\label{lem:ChNorm}
  We have
 \begin{equation*}
  \|\Ch\bv\|_\aP\leq \PFfactor \|\bv\|_\aP\qquad\forall\,\bv\in \DESh.
 \end{equation*}
\end{lemma}
\begin{proof}
 Let $\qz=\Ch\bv$.  Then $\qz$ (and hence $\qmz$) belongs to $\KERP$.
  It follows from \eqref{eq:mCoercive}, \eqref{eq:Bounded} and \eqref{eq:ChDef}  that
\begin{align*}
   \|\Ch\bv\|_\aP^2=\|\qz\|_{a\times a}^2
                           &=\cB(\qz,\qmz)\\
                           &=\cB(\bv,\qmz)
                          \leq\PFfactor\|\bv\|_\aP\|\Ch\bv\|_\aP.
\end{align*}
\end{proof}
\begin{corollary}\label{cor:VHNormEquivalence}
  The following relations are valid$\,:$
\begin{alignat}{3}
 \|\bv-\Ch\bv\|_\aP&\leq [2+(\CPF/\alpha)]\|\bv\|_\aP
 &\qquad&\forall\,\bv\in \DESh,
 \label{eq:VHNormEquivalence1}\\
 \|\bv\|_\aP&\leq C_\dag\sqrt{\beta/\alpha}\|\bv-\Ch\bv\|_\aP
 &\qquad&\forall\,\bv\in \DESH.
 \label{eq:VHNormEquivalence2}
\end{alignat}
\end{corollary}
\begin{proof} The  inequality \eqref{eq:VHNormEquivalence1}
 follows from Lemma~\ref{lem:ChNorm} and the triangle inequality, and
 the inequality \eqref{eq:VHNormEquivalence2} follows from \eqref{eq:Projection},
  Remark~\ref{rem:AlternativePiHEst} and \eqref{eq:Ch1} :
\begin{equation*}
  \|\bv\|_\aP=\|(\PiH\times\PiH) (\bv-\Ch\bv)\|_\aP\leq
   C_\dag\sqrt{\beta/\alpha}\,\|\bv-\Ch\bv||_\aP.
\end{equation*}
\end{proof}
\subsection{The Finite Element Space $\MV$}\label{subsec:MV}
 The ideal multiscale finite element space
  $$\MV\subset \DESh$$
 is defined by
\begin{equation}\label{eq:MVDef}
  \MV=\{\bv\in \DESh:\,\B{\bv}{\bw}=0\quad\forall\,\bw\in \KERP\}.
\end{equation}
\par
 Let $\bv\in\DESh$ be arbitrary.
 It follows from Lemma~\ref{lem:BStability} (with $V=\KER$),
 \eqref{eq:ChDef} and \eqref{eq:MVDef} that
\begin{equation}\label{eq:MVCharacterization}
  \bv\in \MV\quad\Leftrightarrow\quad \cB(\Ch\bv,\bw)=0
  \quad \forall\,\bw\in \KERP
 \quad \Leftrightarrow\quad \Ch\bv=0.
\end{equation}
 Therefore we have $\MV=(\bI-\Ch)(\DESh)$ and
\begin{equation}\label{eq:dimMV}
  \mathrm{dim}\MV=2n-2\ell=2m.
\end{equation}
\par
 A basis for $\MV$ is given in the following lemma.
\begin{lemma}\label{lem:MSBasis}
 Let $\{\phi_1,\ldots,\phi_m\}$ be the nodal basis of $\VH$,
 $\bpsi_i=\Ch(\phi_i,0)$ and $\bxi_i=\Ch(0,\phi_i)$.  Then
  $$\{(\phi_1,0)-\bpsi_1,\ldots,(\phi_m,0)-\bpsi_m,
  (0,\phi_1)-\bxi_1,\ldots,(0,\phi_m)-\bxi_m\}$$
   is a basis for $\MV$.
\end{lemma}
\begin{proof} In view of \eqref{eq:dimMV}, it suffices to show
 that  the $2m$ functions
 $(\phi_i,0)-\bpsi_i=(\bI-\Ch)(\phi_i,0)$ and
 $(0,\phi_i)-\bxi_i=(\bI-\Ch)(0,\phi_i)$
   ($1\leq i\leq m$)
 are linearly independent.
\par
  Suppose   $\sum_{i=1}^m \big[c_i(\bI-\Ch)(\phi_i,0)
  +d_i(\bI-\Ch)(0,\phi_i)\big]=\bf{0}$.
   It then follows from
  \eqref{eq:Projection} and \eqref{eq:Ch1} that $\sum_{i=1}^m \big[c_i(\phi_i,0)
  +d_i(0,\phi_i)\big]=\bf{0}$ and hence
  $c_i=d_i=0$ for $1\leq i\leq m$.
\end{proof}
\begin{remark}\label{rem:MV}\rm
 It follows from Lemma~\ref{lem:MSBasis} that we also have
  $\MV=(\bI-\Ch)(\DESH)$.
\end{remark}
\begin{remark}\label{rem:Symmetry}\rm
  Let $(\psi_{i,1},\psi_{i,2})=\bpsi_i=\Ch(\phi_i,0)$ and
  $(\xi_{i,1},\xi_{i,2})=\bxi_i=\Ch(0,\phi_i)$.
  It follows from \eqref{eq:ChDef} and the relation (cf. \eqref{eq:cBDef})
   $$\B{(y,-p)}{\qz}=\B{\py}{(z,-q)} \qquad\forall\,\py,\qz\in \ES$$
   that
  $\psi_{i,1}=\xi_{i,2}$ and $\psi_{i,2}=-\xi_{i,1}$.
\end{remark}
\begin{figure}[htbp]
\centering
\includegraphics[width=.4\linewidth]{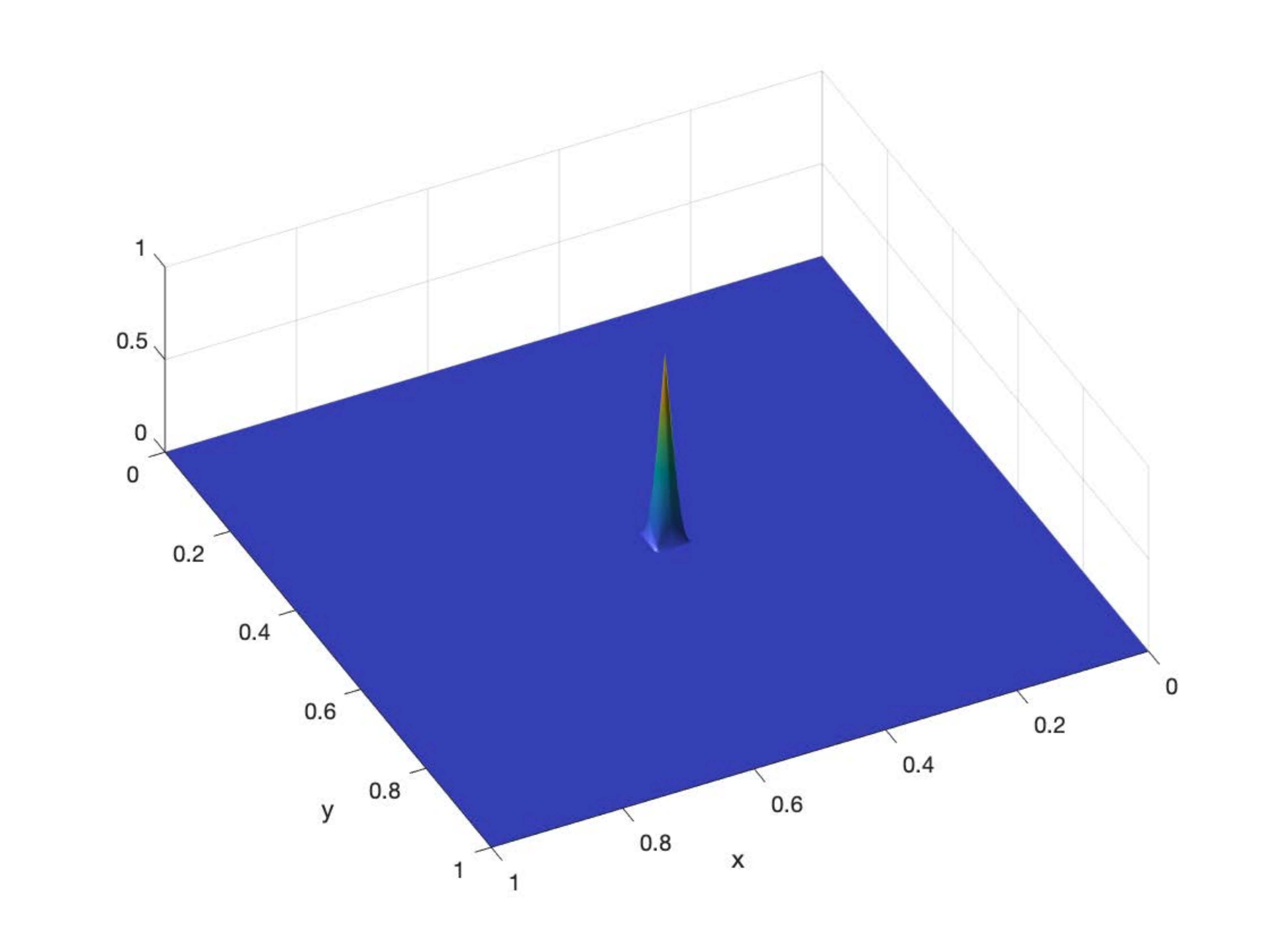}
\caption{A basis function $\phi_i$ of $V_H$.}
\label{figure:StandardBasis}
\end{figure}
\par
\begin{figure}[htbp]
\centering
\includegraphics[width=.4\linewidth]{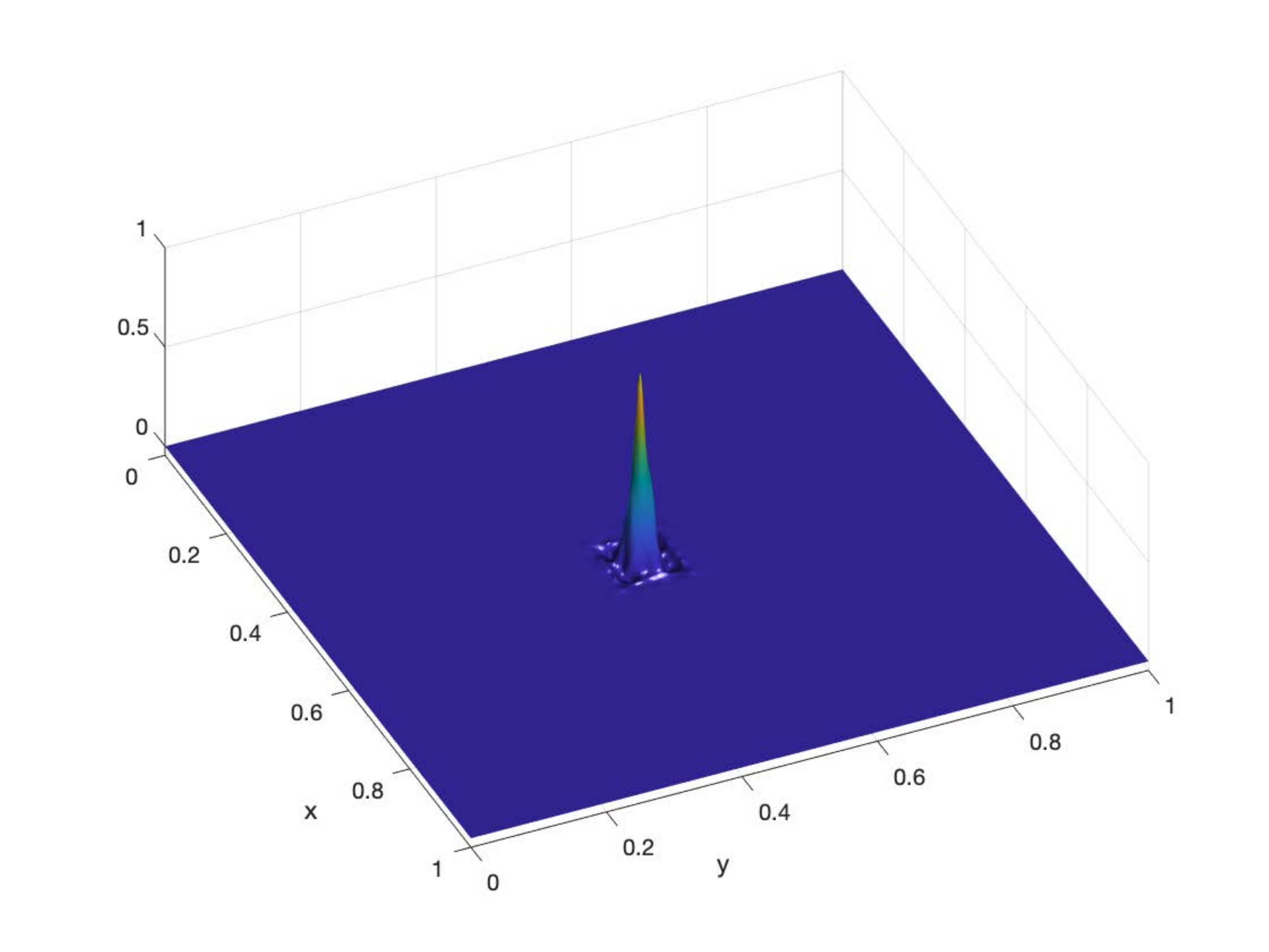}
\quad\includegraphics[width=.4\linewidth]{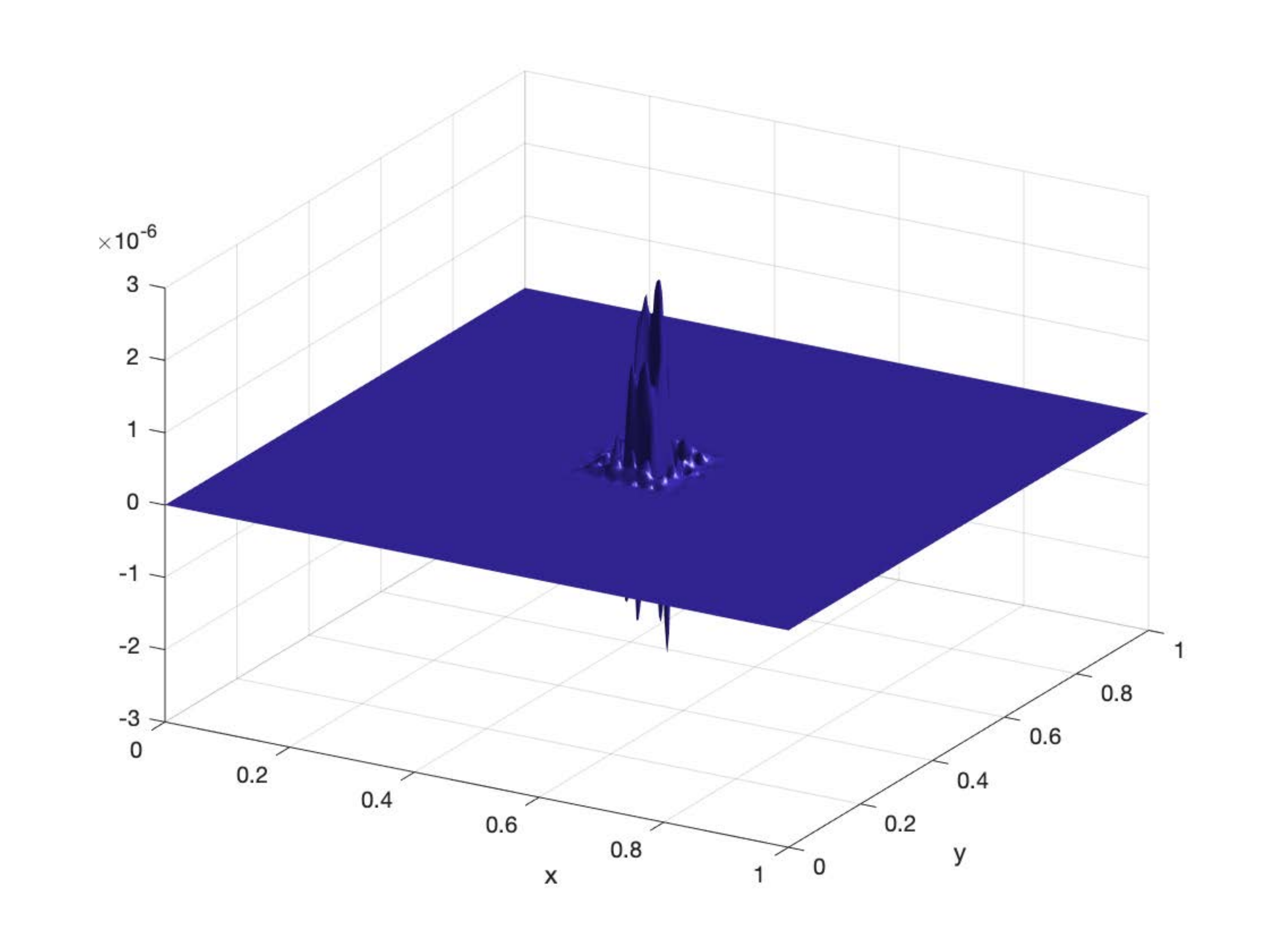}
\caption{The basis function $(\phi_i,0)-\Ch(\phi_i,0)$ of $\MV$ for
 Example~\ref{example:Oscillatory} with $H=1/40$, $h=1/320$ and
 $\epsilon=0.0025$: first component (left) and
second component (right)}
\label{figure:PsiOscillatory}
\end{figure}
\begin{figure}[htbp]
\centering
\includegraphics[width=.4\linewidth]{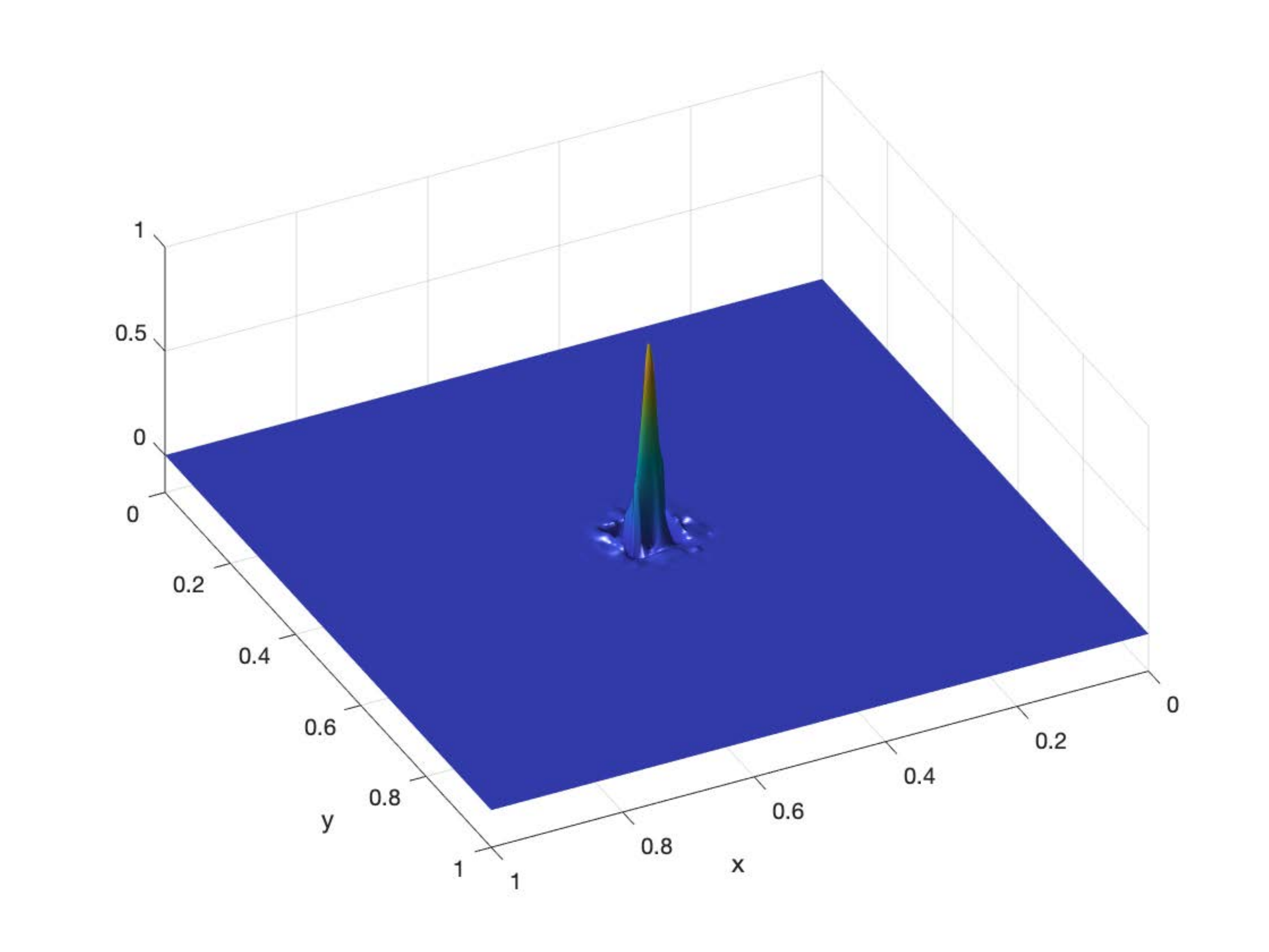}
\quad\includegraphics[width=.4\linewidth]{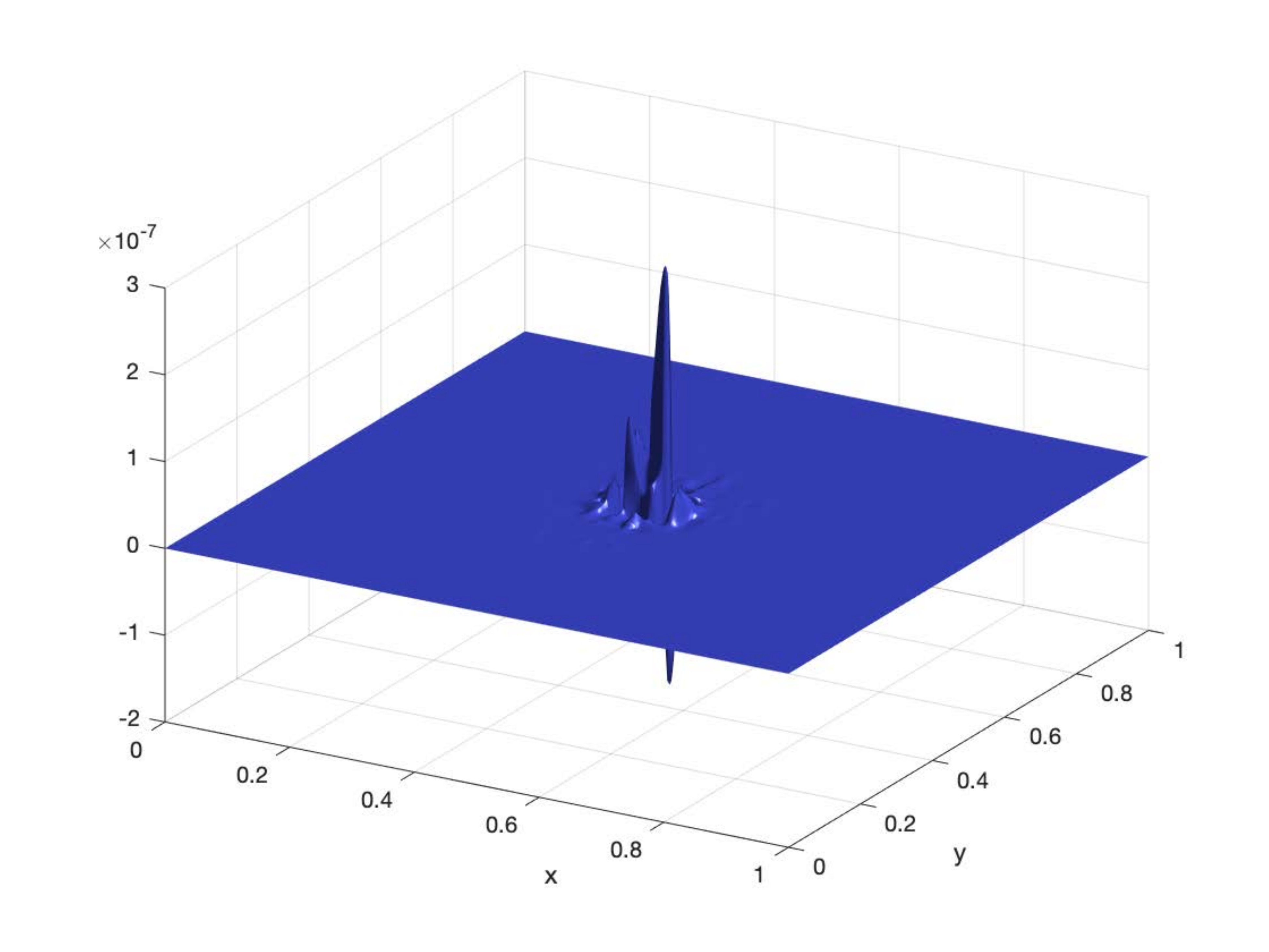}
\caption{The basis function $(\phi_i,0)-\Ch(\phi_i,0)$ of $\MV$ for
 Example~\ref{example:Heterogeneous} with $H=1/40$ and $h=1/320$: first component (left) and
second component (right)}
\label{figure:PsiHeterogeneous}
\end{figure}

 The figure of a typical basis function $\phi_i$ of $V_H$ is given  in
 Figure~\ref{figure:StandardBasis}, and the figure of the corresponding
 basis function $(\phi_i,0)-\Ch(\phi_i,0)$
 is displayed in Figure~\ref{figure:PsiOscillatory} for Example~\ref{example:Oscillatory}, and
 in Figure~\ref{figure:PsiHeterogeneous}
 for Example~\ref{example:Heterogeneous}.
\begin{remark}\label{rem:ED}\rm
 The exponential decay of $\bpsi_i=\Ch(\phi_i,0)$ (and hence $\bxi_i=\Ch(0,\phi_i)$
 in view of Remark~\ref{rem:Symmetry}) are
 clearly observed in
 Figure~\ref{figure:PsiOscillatory} and Figure~\ref{figure:PsiHeterogeneous}.
\end{remark}
%
\subsection{The Discrete Problem}
 The approximate solution $\MSpy\in\MV$ is defined by
\begin{equation}\label{eq:IdealMS}
  \B{\MSpy}{\qz}=\int_\O y_dq\,dx\qquad\forall\,\qz\in\MV.
\end{equation}
\par
  The well-posedness of \eqref{eq:IdealMS} is guaranteed by the following lemma.
\begin{lemma}\label{lem:Stability}
  We have
\begin{align}
  \inf_{\bv\in \MV}\sup_{\bw\in\MV}\frac{\B{\bv}{\bw}}
  {\|\bv\|_\aP\|\bw\|_\aP}&\geq [2+(\CPF/\alpha)]^{-1}
    \label{eq:Stability}.
\end{align}
\end{lemma}
\begin{proof}
  Let $\bv=\qz\in\MV$ be arbitrary.  Then $\qmz-\Ch\qmz\in \MV$
  and
 it follows from \eqref{eq:mCoercive},
 \eqref{eq:VHNormEquivalence1} and \eqref{eq:MVDef}
 that
\begin{align*}
  \|\bv\|_\aP &=\frac{\B{\qz}{\qmz}}{\|\qmz\|_\aP}\\
  &=\frac{\B{\qz}{\qmz-\Ch\qmz}}{\|\qmz\|_\aP}\\
  &\leq
    [2+(\CPF/\alpha)]
  \bigg(\frac{\B{\qz}{\qmz-\Ch\qmz}}{\|\qmz-\Ch\qmz\|_\aP}\bigg)\\
  &\leq [2+(\CPF/\alpha)]\sup_{\bw\in\MV}\frac{\B{\bv}{\bw}}{\|\bw\|_\aP}
\end{align*}
\end{proof}
\subsection{Energy Error}
 It follows from \eqref{eq:FineDiscreteProblem} and \eqref{eq:IdealMS} that
\begin{equation}\label{eq:Galerkin}
  \B{\pyh-\MSpy}{\qz}=0\qquad\forall\,\qz\in \MV.
\end{equation}
 We will use the Galerkin relation \eqref{eq:Galerkin} to derive an error estimate for
 the ideal multiscale finite element method defined by \eqref{eq:IdealMS}.
\begin{theorem}\label{thm:IdealEnergyError}
  We have
\begin{equation}\label{eq:pyEnergyError}
  \|\pyh-\MSpy\|_\aP\leq (C_\dag/\sqrt{\alpha})H\|y_d\|_\LT.
\end{equation}
\end{theorem}
\begin{proof}
 In view of Remark~\ref{rem:Isomorphism}, \eqref{eq:MVDef} and \eqref{eq:Galerkin},
  we have $ \pyh-\MSpy\in\KERP$ and consequently
\begin{equation}\label{eq:MSSolution}
  p_h-\MSp \quad \text{and} \quad y_h-\MSy \quad\text{belong to} \quad \KER.
\end{equation}
\par
 Putting \eqref{eq:Comparison}, \eqref{eq:mCoercive},
 \eqref{eq:FineDiscreteProblem}, \eqref{eq:PiEstimates},
 \eqref{eq:MVDef} and
 \eqref{eq:MSSolution} together, we have
\begin{align*}
 \|\pyh-\MSpy\|_\aP^2
  &=\B{\pyh-\MSpy}{\pmyh-\MSpmy}\\
    &=\B{\pyh}{\pmyh-\MSpmy}\\
  &=\int_\O y_d(p_h-\MSp)dx\\
  &=\int_\O y_d[(p_h-\MSp)-\PiH(p_h-\MSp)]dx\\
  &\leq C_\dag H\|y_d\|_\LT |p_h-\MSp|_\HO\\
  &\leq (C_\dag/\sqrt{\alpha})H\|y_d\|_\LT\|p_h-\MSp\|_a,
\end{align*}
 and \eqref{eq:pyEnergyError} follows immediately.
\end{proof}
\begin{remark}\label{rem:ErrorRepresentation}\rm
  In view of \eqref{eq:Ch2}, \eqref{eq:MVCharacterization} and
   \eqref{eq:MSSolution}, we can express the error
  of the ideal multiscale finite element method as
\begin{equation*}
  \pyh-\MSpy=\Ch(\pyh-\MSpy)=\Ch\pyh.
\end{equation*}
\end{remark}
%
\subsection{$L_2$ Error}
 We will obtain an estimate for the $L_2$ error by a duality argument.
\begin{theorem}\label{thm:IdealLTwoError}
  We have
\begin{equation}\label{eq:IdealLTError}
  \|\pyh-\MSpy\|_{\LT\times\LT}\leq  \PFfactor(C_\dag^2/{\alpha})H^2\|y_d\|_\LT.
\end{equation}
\end{theorem}
\begin{proof}
 In view of Remark~\ref{rem:Isomorphism}, we can define
 $\qz\in \DESh$ by
\begin{equation}\label{eq:Duality1}
 \B{\qz}{\rs}=\int_\O (p_h-\MSp)r\,dx+\int_\O (y_h-\MSy)s\,dx\qquad\forall\,
 \rs\in  \DESh.
\end{equation}
\par
 Let  $\Ch\qz\in\KERP$ be denoted by $\tqz$.  It follows from \eqref{eq:Comparison},
 \eqref{eq:mCoercive}, \eqref{eq:PiEstimates},
 \eqref{eq:ChDef}, \eqref{eq:Ch5}, \eqref{eq:Ch2},  \eqref{eq:Duality1}
 and the Cauchy-Schwarz inequality that
\begin{align*}
  \|\tqz\|_\aP^2&=\B{\tqz}{\tqmz}\\
  &=\B{\qz}{\tqmz}\\
  &=\int_\O (p_h-\MSp)\tq\,dx-\int_\O (y_h-\MSy)\tz\,dx\\
  &= \int_\O(p_h-\MSp)(\tq-\PiH\tq)dx-\int_\O (y_h-\MSy)(\tz-\PiH\tz)dx\\
  &\leq C_\dag H\big(\|p_h-\MSp\|_\LT |\tq|_{H^1(\O)}
  +\|y_h-\MSy\|_\LT|\tz|_{H^1(\O)}\big)\\
  &\leq C_\dag H \|(p_h-\MSp,y_h-\MSy)\|_{\LT\times\LT}
  (1/\sqrt{\alpha})\|\tqz\|_\aP,
\end{align*}
 and hence
\begin{equation}\label{eq:Duality2}
   \|\Ch\qz\|_\aP
    \leq (C_\dag/\sqrt\alpha)H \|\pyh-\MSpy\|_{\LT\times\LT}.
\end{equation}
\par
 On the other hand, we have
\begin{align}\label{eq:Duality3}
  &\|(p_h-\MSp,y_h-\MSy)\|_{\LT\times\LT}^2=\B{\qz}{(p_h-\MSp,y_h-\MSy)}\notag\\
      &\hspace{40pt}=\B{\Ch\qz}{(p_h-\MSp,y_h-\MSy)}\\
      &\hspace{40pt}\leq\PFfactor\|\Ch\qz\|_\aP\|(p_h-\MSp,y_h-\MSy)\|_\aP\notag
\end{align}
 by \eqref{eq:Bounded}, \eqref{eq:ChDef}, \eqref{eq:MSSolution}
 and \eqref{eq:Duality1}.
\par
 Putting \eqref{eq:pyEnergyError}, \eqref{eq:Duality2} and \eqref{eq:Duality3} together, we
 arrive at
\begin{align*}
 &\|\pyh-\MSpy\|_{\LT\times\LT}\notag\\
 &\hspace{60pt}\leq \PFfactor(C_\dag/\sqrt{\alpha})H
 \|(p_h-\MSp,y_h-\MSy)\|_\aP\\
   &\hspace{60pt}\leq  \PFfactor(C_\dag^2/{\alpha})H^2\|y_d\|_\LT.\notag
\end{align*}
\end{proof}
\section{A Localized Multiscale Finite Element Space}\label{sec:LMFES}
 The constructions of $\bpsi_i=\Ch(\phi_i,0)\in\KERP$ and
 $\bxi_i=\Ch(0,\phi_i)\in \KERP$ require solving the equations
\begin{alignat}{3}
  \B{\bpsi_i}{\qz}&=\B{(\phi_i,0)}{\qz}\qquad\forall\,
  \qz\in\KERP,\label{eq:Corrector1}\\
   \B{\bxi_i}{\qz}&=\B{(0,\phi_i)}{\qz}\qquad\forall\,\qz\in\KERP,
   \label{eq:Corrector2}
\end{alignat}
 which are expensive.  However the exponential decays of $\bpsi_i$ and $\bxi_i$ observed in
  Figures~\ref{figure:PsiHeterogeneous}, Figures~\ref{figure:PsiOscillatory} and Remark~\ref{rem:ED}
 indicate that  it is possible to capture $\bpsi_i$ and $\bxi_i$ by local approximations.
 (Note that in practice we only need to solve one of these equations because of the observation in
 Remark~\ref{rem:Symmetry}.)
\par
 We will construct a localized multiscale
 finite element space
 by replacing $\bpsi_i$ (resp., $\bxi_i$) with an approximate solution
 of \eqref{eq:Corrector1} (resp.,  \eqref{eq:Corrector2})
 obtained by a preconditioned minimum residual (P-MINRES) algorithm
 (cf. \cite[Chapter~8]{Greenbaum:1997:Iterative} and
 \cite[Section~4.1]{ESW:2014:FEFIS}).
 Our construction extends those in \cite{KPY:2018:LOD,BGS:2021:LOD}
 to  symmetric indefinite problems.
\subsection{An Additive Schwarz Preconditioner}\label{subeq:ASP}
 Let $\AKER:\KER\longrightarrow (\KER)'$ be the linear operator defined by
\begin{equation}\label{eq:AKERDef}
  \langle \AKER v,w\rangle=a(v,w) \qquad\forall\,v,w\in\KER,
\end{equation}
 where $a(\cdot,\cdot)$ is given in \eqref{eq:aDef}.
 We begin by constructing an additive Schwarz preconditioner (cf. \cite{DW:1987:AS,TW:2005:DD,BScott:2008:FEM})
  for $\AKER$.
\par
\begin{figure}[ht]
\begin{center}
  \includegraphics[width=3in]{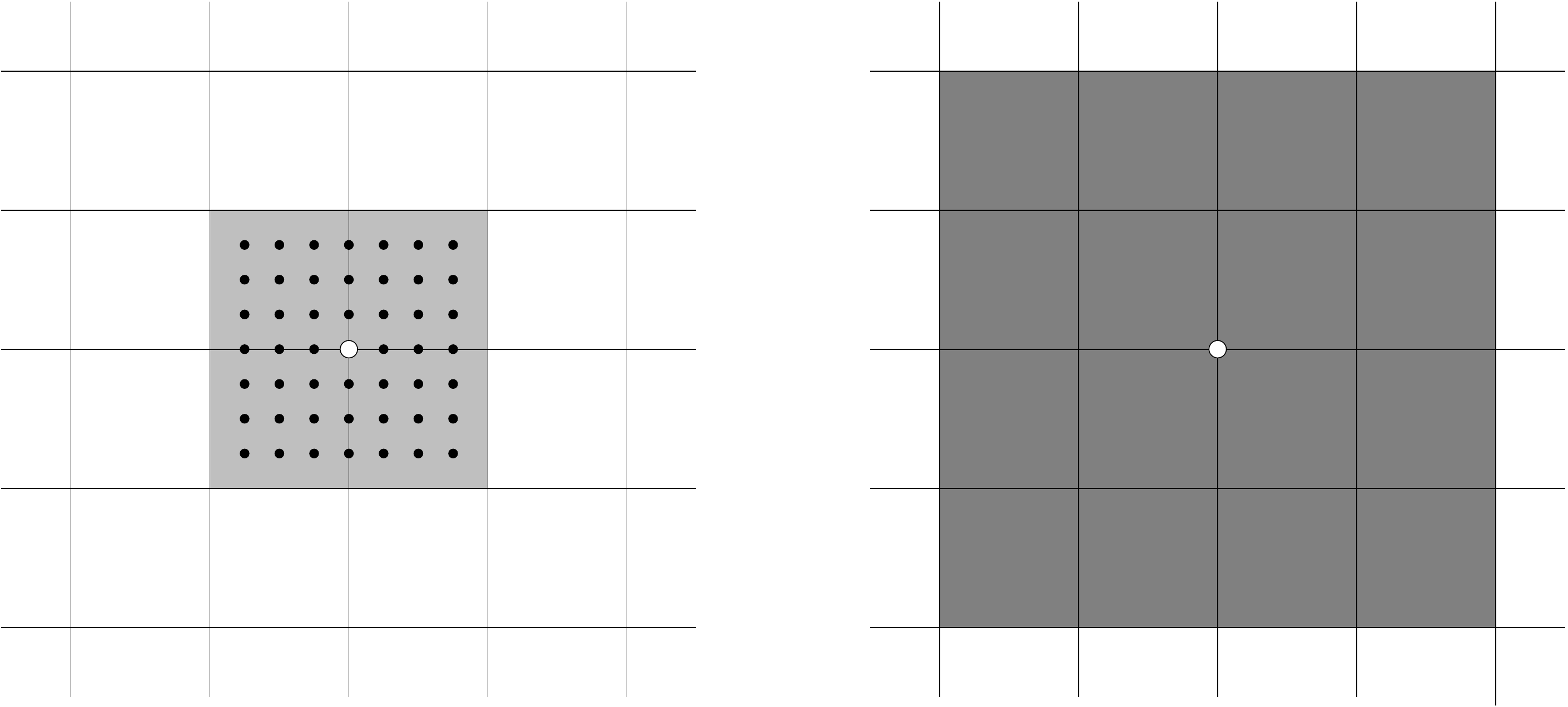}
\end{center}
\caption{The patches $\omega_{x_i}$ $($left$)$
and $\tilde{\omega}_{x_i}$ $($right$)$,
the node $x_i$ is represented by the
 circle and the nodes for $\varphi_{i,1},\ldots,\varphi_{i,m_i}$
 are represented by the
 solid dots. }
\label{fig:Patch}
\end{figure}

 Let $x_1,\ldots,x_m$ be the (interior) nodes for $V_H$.  We define the
 subspaces $\Ki$ ($1\leq i\leq m$)
 of $\KER$ by
\begin{equation}\label{eq:KiDef}
  \Ki=\{(I-\PiH)v:\,\text{$v\in V_h$ and $v$ vanishes outside $\omega_{x_i}$}\},
\end{equation}
 where $\omega_{x_i}$ is the union of the elements in $\cT_H$ that share $x_i$
 as a common vertex
 (cf. Figure~\ref{fig:Patch} for a two dimensional example with the $Q_1$ element).
 The functions in $\Ki$ are supported
 on the patch $\tilde\omega_{x_i}$ obtained from $\omega_{x_i}$ by
 adding one layer of elements in
 $\cT_H$ (cf. Figure~\ref{fig:Patch}).  Let $\varphi_{i,1},\ldots,\varphi_{i,m_i}$
 be the nodal basis
 functions of $V_h$ that vanish at $x_i$ and outside $\omega_{x_i}$
  (cf. Figure~\ref{fig:Patch}).  Then, as in
 Lemma~\ref{lem:KERBasis},
 $\{(I-\PiH)\varphi_{i,1},\ldots,(I-\PiH)\varphi_{i,m_i}\}$ is a basis of $\Ki$.
\par
 Let $I_i:\Ki\longrightarrow \KER$ be the natural injection.
 The SPD additive Schwarz preconditioner
  $\SKER:(\KER)'\longrightarrow \KER$ for $\AKER$ is given by
\begin{equation}\label{eq:SDef}
 \SKER=\sum_{i=1}^m I_i (\AKERi)^{-1} I_i^t,
\end{equation}
 where $\AKERi:\Ki\longrightarrow (\Ki)'$ is defined by
\begin{equation}\label{eq:AiDef}
  \langle \AKERi v,w\rangle=a(v,w)\qquad\forall\,v,w\in \Ki.
\end{equation}
\par
 According to the Raleigh quotient formulas, we have
\begin{alignat}{3}
 a(v,v)&=\langle \AKER v,v\rangle\leq \LMax \langle(\SKER)^{-1}v,v\rangle
 &\qquad& \forall\,v\in \KER, \label{eq:RMax}\\
 a(v,v)&=\langle \AKER v,v\rangle\geq \LMin \langle(\SKER)^{-1}v,v\rangle
 &\qquad& \forall\,v\in \KER, \label{eq:RMin}
\end{alignat}
 and the  following spectral estimates
 can be found in \cite[Section~3]{BGS:2021:LOD}:
\begin{equation}\label{eq:SpectarlEstimates}
  \LMax\leq C_{\rm upper} \quad\text{and}\quad
   \LMin\geq C_{\rm lower}(\alpha/\beta),
\end{equation}
 where the positive constants $C_{\rm upper}$ and $C_{\rm lower}$
 only depend on the shape regularity of $\cT_H$.
\subsection{The Generalized Finite Element Space $\MSV$}\label{subsec:MSV}
 Let $$\bB:\KERP\longrightarrow (\KERP)'$$ be the linear operator defined by
\begin{equation}\label{eq:bBDef}
  \langle \bB \bv,\bw\rangle=\cB(\bv,\bw)\qquad\forall\,\bv,\bw\in\KERP.
\end{equation}
 We can then rewrite \eqref{eq:Corrector1} and \eqref{eq:Corrector2} as
\begin{align}
  \bB\bpsi_i&=\bm{f}_i,\label{eq:Corrector12}\\
  \bB\bxi_i&=\bm{g}_i,\label{eq:Corrector22}
\end{align}
 where $\bm{f}_i,\bm{g}_i\in (\KERP)'$ are defined by
\begin{equation*}
  \langle\bm{f}_i,\bw\rangle=\B{(\phi_i,0)}{\bw}
  \quad\text{and}\quad
   \langle\bm{g}_i,\bw\rangle=\B{(0,\phi)}{\bw}\qquad\forall\,\bw\in \KERP.
\end{equation*}
\par
 Let $\bpsi_{i,k}\in\KERP$ (resp., $\bxi_{i,k}\in\KERP$) be the output of
 $k$ steps of the P-MINRES algorithm
 applied  to
 \eqref{eq:Corrector12} (resp., \eqref{eq:Corrector22}) with initial
 guess $0$, where the SPD preconditioner
 $$\bS:(\KERP)'\longrightarrow \KERP$$ is given by
\begin{equation}\label{eq:bSDef}
  \bS(\mu,\rho)=(\SKER\mu,\SKER\rho). 
\end{equation}
 Then the $2m$  functions
\begin{equation}\label{eq:LocalizedBasis}
 (\phi_1,0)-\bpsi_{1,k},\ldots,(\phi_m,0)-\bpsi_{m,k},(0,\phi_1)-\bxi_{1,k},
 \ldots,(0,\phi_m)-\bxi_{m,k}
\end{equation}
 are linearly independent because the intersection of $\VH\times \VH$ and
  $\KERP$ is trivial, and we
 define the  generalized finite element space $\MSV$ by
\begin{align}\label{eq:MSVDef}
  \MSV&=\mathrm{span}\big\{ (\phi_1,0)-\bpsi_{1,k},\ldots,(\phi_m,0)-\bpsi_{m,k},
  (0,\phi_1)-\bxi_{1,k},\ldots,(0,\phi_m)-\bxi_{m,k}\big\}.
\end{align}
\begin{remark}\label{rem:Localization}\rm
 It follows from \eqref{eq:cBDef},
 \eqref{eq:SDef} and \eqref{eq:bBDef} that the support of $\bS\bB\qz$ is a subset
 of the union of all the $\tilde\omega_{x_i}$
 whose intersection with the support of $\qz$ have nonempty interiors.
 As the output of a preconditioned
 Krylov subspace method
 with initial guess $0$, the function $\bpsi_{i,k}$ belongs to
    $$\mathrm{span}\big\{\bS (\phi_i,0), (\bS \bB)\bS(\phi_i,0),
    \ldots,(\bS\bB)^k\bS(\phi_i,0)\big\},    $$
 and hence is supported in a patch around $x_i$ (with respect to
 $\cT_H$) whose diameter is proportional to $k H$.  This is also true for the function
 $\bxi_{i,k}$.
 The functions in \eqref{eq:LocalizedBasis} are therefore
 locally corrected basis functions and $\MSV$ defined by \eqref{eq:MSVDef}
  is a localized
 multiscale finite element space.
\end{remark}
\begin{remark}\label{rem:bBbS} \rm
 In view of \eqref{eq:cBDef}, we can express $\bB$ in the matrix form
\begin{equation}\label{eq:bBMatrix}
  \bB=\begin{bmatrix}
    \AKER &\MKER \\[4pt]
    \MKER&-\AKER
  \end{bmatrix},
\end{equation}
 where $\MKER:\KER\longrightarrow (\KER)'$ is the (symmetric) linear
 operator defined by
\begin{equation*}
  \langle\MKER r,z\rangle=\int_\O rz\,dx \qquad\forall\,r,z\in \KER.
\end{equation*}
 We can also write $\bS$ as the diagonal matrix
\begin{equation}\label{eq:bSMatrix}
  \bS=\begin{bmatrix}
    \SKER & 0\\[4pt]
    0&\SKER
  \end{bmatrix}.
\end{equation}
\end{remark}
%
\subsection{Spectral Analysis of $\bS\bB:\KERP\longrightarrow\KERP$}
 A spectral analysis of the operator $\bS\bB$ is provided
  in the following lemma, which is the key to the analysis of the
 localized multiscale finite element method in Section~\ref{sec:LMFEM}.
\begin{lemma}\label{lem:SpectralAnalysis}
  The spectrum $\sigma(\bS\bB)$ of $\,\bS\bB$ satisfies
\begin{equation}\label{eq:Spectrum}
  \sigma(\bS\bB)\subseteq [-d_*,-c_*]\cup[c_*,d_*],
\end{equation}
 where
\begin{equation}\label{eq:Constants}
 \text{\rm $c_*= \LMin$ and
 $d_*=\LMax[1+(\CPF/\alpha)]$.}
\end{equation}
\end{lemma}
\begin{proof}
   In view of \eqref{eq:bBMatrix} and \eqref{eq:bSMatrix},
   the eigenvalues of $\bS\bB$ are real numbers.
   Let $\lambda$ be one of the eigenvalues and $\rs\in\KERP$ be
   a corresponding eigenvector.
   Given any $\qz\in\KERP$, we have, by \eqref{eq:bBDef}, 
\begin{equation}\label{eq:EigenVector}
 \cB(\rs,\qz)=\langle (\bS)^{-1}\bS\bB\rs,\qz\rangle
                  =\lambda\langle(\bS)^{-1}\rs,\qz\rangle.
\end{equation}
\par
 It follows from \eqref{eq:mCoercive},
 \eqref{eq:AKERDef}, \eqref{eq:RMin} and
 \eqref{eq:EigenVector} that
\begin{align*}
 \|\rs\|_{a\times a}^2=\cB(\rs,\rms)
    &=\lambda\langle(\bS)^{-1}\rs,\rms\rangle\\
    &=\lambda\big(\langle(\SKER)^{-1}r,r\rangle-
    \langle(\SKER)^{-1}s,s\rangle\big)\\
    &\leq |\lambda|\big(\langle(\SKER)^{-1}r,r\rangle+
    \langle(\SKER)^{-1}s,s\rangle\big)\\
     &=|\lambda|[1/\LMin]\|\rs\|_{a\times a}^2
\end{align*}
 and hence
\begin{equation*}
 |\lambda|\geq \LMin.
\end{equation*}
\par
 On the other hand we can deduce from \eqref{eq:Bounded},
 \eqref{eq:RMax} and \eqref{eq:EigenVector}
 that
\begin{align*}
| \lambda| \|\rs\|_{a\times a}^2
&\leq \LMax \big|\lambda \big(\langle (\SKER)^{-1} r,r\rangle
+\langle (\SKER)^{-1} s,s\rangle\big)\big|\\
     &= \LMax \big|\lambda\langle(\bS)^{-1}\rs,\rs\rangle\big|\\
     &= \LMax |\cB(\rs,\rs)|\\
     &\leq \LMax [1+(\CPF/\alpha)]\|\rs\|_\aP^2
\end{align*}
 and hence
\begin{equation*}
 |\lambda|\leq \LMax\PFfactor .
\end{equation*}
\end{proof}
\begin{corollary}\label{cor:NormComparison}
   The following relations hold for any $\bv\in\KERP :$
\begin{equation}\label{eq:NormComparison}
  \bigg(\frac{c_*}{\LMax^\frac12}\bigg)\|\bv\|_\aP\leq
  \langle\bB \bv,\bS\bB\bv\rangle^\frac12\leq
  \bigg(\frac{d_*}{\LMin^\frac12}\bigg)\|\bv\|_\aP.
\end{equation}
\end{corollary}
\begin{proof}
  Since the operator $\bS\bB$ is symmetric with respect to the
  inner product
  defined by
  $\langle (\bS)^{-1}\cdot,\cdot\rangle$, we have, by
  Lemma~\ref{lem:SpectralAnalysis} and the  spectral theorem,
\begin{equation*}
  c_* \langle(\bS)^{-1}\bv,\bv\rangle^\frac12\leq
  \langle\bB \bv,\bS\bB\bv\rangle^\frac12
  \leq d_*\langle (\bS)^{-1}\bv,\bv\rangle^\frac12,
\end{equation*}
 which together with
 \eqref{eq:RMax} and \eqref{eq:RMin} implies \eqref{eq:NormComparison}.
\end{proof}
\section{The Localized Multiscale Finite Element Method} \label{sec:LMFEM}
 The localized multiscale finite element method is to find
 $\LMSpy\in\MSV$ such that
\begin{equation}\label{eq:LMFEM}
 \cB\big(\LMSpy,\qz\big)=\int_\O y_d q\,dx\qquad\forall\,\qz\in\MSV.
\end{equation}
\par
 We will keep track of all the constants in the error analysis so that the
 constants that appear in the energy error estimate (cf. Theorem~\ref{thm:LMFEMEnergyError}) and the
 $\LTP$ error estimate (cf. Theorem~\ref{thm:LMFEMLTwoError}) are independent of the mesh sizes
 ($h$ and $H$) and the contrast
 $\beta/\alpha$.
\par
 We begin the error analysis by comparing $\bpsi_i$
 and $\bpsi_{i,k}$.
 (resp., $\bxi_i$ and $\bxi_{i,k}$).
%
\subsection{The Relation between $(\bpsi_i,\bxi_i)$ and $(\bpsi_{i,k},\bxi_{i,k}$)}
\label{subsec:Relations1}
\par
 It follows from Lemma~\ref{lem:SpectralAnalysis} and the theory of the P-MINRES
 algorithm (cf. \cite[Theorem~4.14]{ESW:2014:FEFIS}) that
\begin{align}
 \langle \bB(\bpsi_i-\bpsi_{i,k}),\bS\bB(\bpsi_i-\bpsi_{i,k})\rangle^\frac12&\leq
   2 q^{\lfloor k/2\rfloor}
   \langle \bB\bpsi_{i},\bS\bB\bpsi_{i}\rangle^\frac12,\label{eq:ResidualError1}\\
    \langle \bB(\bxi_i-\bxi_{i,k}),\bS\bB(\bxi_i-\bxi_{i,k})\rangle^\frac12&\leq
   2q^{\lfloor k/2\rfloor}
   \langle \bB\bxi_{i},\bS\bB\bxi_{i}\rangle^\frac12,\label{eq:ResidualError2}
\end{align}
 where (cf. \eqref{eq:Constants})
\begin{equation}\label{eq:qDef}
 q=\frac{d_*-c_*}{d_*+c_*}\\
    =\frac{\kappa(\SKER\AKER)\PFfactor-1}{\kappa(\SKER\AKER)\PFfactor+1}.
\end{equation}
\begin{remark}\label{rem:ConditionNUmber}\rm
  It follows from \eqref{eq:SpectarlEstimates} that the condition number
  $\kappa(\SKER\AKER)$ has the following (pessimistic)
  upper bound that is independent of the mesh sizes:
\begin{equation}\label{eq:ConditionNumber}
   \Cond=\frac{\LMax}{\LMin}
   \leq \Big(\frac{C_{\rm upper}}{C_{\rm lower}}\Big)\Big(\frac{\beta}{\alpha}\Big).
\end{equation}
 Consequently we also have the bound
\begin{equation}\label{eq:qBound}
  q\leq \frac{\d\Big(\frac{C_{\rm upper}}{C_{\rm lower}}\Big)\Big(\frac{\beta}{\alpha}\Big)\PFfactor-1}
  {\d\Big(\frac{C_{\rm upper}}{C_{\rm lower}}\Big)\Big(\frac{\beta}{\alpha}\Big)\PFfactor+1}
\end{equation}
 that is independent of the mesh sizes, but which may be too pessimistic.
\end{remark}
\begin{lemma}\label{lem:PMRESError}
  We have
\begin{alignat}{3}
  \|\bpsi_i-\bpsi_{i,k}\|_\aP&\leq C_\star q^{\lfloor k/2\rfloor}\|\bpsi_i\|_\aP
  &\qquad&\text{for}\quad 1\leq i\leq m \label{eq:PMRESError1}\\
\intertext{and}
  \|\bxi_{i}-\bxi_{i,k}\|_\aP&\leq C_\star q^{\lfloor k/2\rfloor}\|\bxi_i\|_\aP
  &\qquad&\text{for}\quad 1\leq i\leq m, \label{eq:PMRESSrror2}
\end{alignat}
 where the positive constant $C_\star$ is  independent of the mesh sizes.
\end{lemma}
\begin{proof} It follows from \eqref{eq:NormComparison} and
\eqref{eq:ResidualError1} that
\begin{align*}
  \|\bpsi_i-\bpsi_{i,k}\|_{a\times a}&\leq
   \frac{\LMax^\frac12}{c_*}
  \langle \bB(\bpsi_i-\bpsi_{i,k}),\bS\bB(\bpsi_i-\bpsi_{i,k})\rangle^\frac12\\
  &\leq 2q^{\lfloor k/2\rfloor}\frac{\LMax^\frac12}{c_*}
  \langle \bB\bpsi_{i},\bS\bB\bpsi_{i}\rangle^\frac12\\
  &\leq 2q^{\lfloor k/2\rfloor}\frac{\LMax^\frac12}{c_*}
  \frac{d_*}{\LMin^\frac12}
    \|\bpsi\|_\aP,
\end{align*}
 which implies \eqref{eq:PMRESError1} with (cf. \eqref{eq:Constants})
\begin{equation}\label{eq:C*Def}
 C_\star =2\Cond^\frac32 \PFfactor.
\end{equation}
\par
 Similarly we can derive \eqref{eq:PMRESSrror2} from \eqref{eq:NormComparison} and \eqref{eq:ResidualError2}.
\end{proof}
\begin{remark}\label{rem:ExponentialDecay}\rm
  It follows from Lemma~\ref{lem:PMRESError} that the basis
   corrections $\Ch(\phi_i,0)=\bpsi_i$ and
  $\Ch(0,\phi_i)=\bxi_i$ ($1\leq i\leq m$) used in the construction of $\MV$
  decay exponentially, as observed in Figure~\ref{figure:PsiHeterogeneous},
  Figure~\ref{figure:PsiOscillatory} and Remark~\ref{rem:ED}.
\end{remark}
\par
 Let the linear operator $\Chk:\DESH\longrightarrow \KERP$ be defined by
\begin{equation}\label{eq:ChkDef}
  \Chk(\phi_i,0)=\bpsi_{i,k} \quad \text{and} \quad \Chk(0,\phi_i)=\bxi_{i,k}
  \quad\text{for}\quad1\leq i\leq m,
\end{equation}
 where $\phi_1,\ldots,\phi_m$ are the nodal basis functions of $V_H$.
\par
 Our next goal is to understand the relation between the operators $\Ch$ and $\Chk$.
%
\subsection{The Relation between $\Ch$ and $\Chk$}
\label{subsec:Relations2}
 We begin with the following lemma.
\begin{lemma}\label{lem:Observations}
 There exist  positive constants
 $C_\diamondsuit$ and $C_\heartsuit$ depending only on the
  shape regularity of $\cT_H$ such that
\begin{equation}\label{eq:Observation1}
  |\phi_i|_\HO\leq C_\diamondsuit H^{\tau_d} \quad\text{with}\quad
    \tau_d=\begin{cases}
      -\frac12 &\qquad d=1\\[2pt]
      \hfill 0&\qquad d=2\\[2pt]
      \hfill\frac12&\qquad d=3
    \end{cases},
\end{equation}
 where $d$ is the dimension of $\O$,
 and
\begin{equation}\label{eq:Observation2}
  \sum_{i=1}^m (|c_i|+|d_i|)\leq C_\heartsuit (1/\sqrt{\alpha})H^{-d}
  \Big\|\sum_{i=1}^m \big[c_i(\phi,0)+d_i(0,\phi_i)\big]\Big\|_\aP
\end{equation}
 for any real numbers $c_1,\ldots,c_m,d_1,\ldots, d_m$,
\end{lemma}
\begin{proof}
  The estimate \eqref{eq:Observation1} follows from a scaling argument.
  For the estimate
  \eqref{eq:Observation2}, we begin with an inverse estimate
\begin{equation*}
  \sum_{i=1}^m(|c_i|+|d_i|)\leq C_\natural H^{-d}
  \Big\|\sum_{i=1}^m \big[c_i(\phi_i,0)+d_i(0,\phi_i)\big]\Big\|_{L_1(\O)\times L_1(\O)},
\end{equation*}
 where the positive constant $C_\natural$ depends only on the shape
 regularity of $\cT_H$.  The proof is then completed
 by the Poincar\'e-Friedrichs inequality
\begin{equation*}
  \|\zeta\|_{L_1(\O)} \leq C_\O |\zeta|_\HO \qquad\forall\,\zeta\in\zHO
\end{equation*}
 and the estimate \eqref{eq:Comparison}.

\end{proof}
\begin{lemma}\label{lem:OperatorDifference}
  The following estimate is valid for any $\bv\in\DESH :$
\begin{equation}\label{eq:OperatorDifference}
  \|(\Ch-\Chk)\bv\|_\aP\leq C_\spadesuit \kappa(\SKER\AKER)^\frac32
   \Big(\frac{\beta}{\alpha}\Big)^\frac12
  \PFfactor^2 H^{-d+\tau_d}q^{\lfloor k/2\rfloor}\|\bv\|_\aP,
\end{equation}
 where the positive constant $C_\spadesuit$ depends only on the shape
 regularity of $\cT_H$.

\end{lemma}
\begin{proof}
   Let $\sum_{i=1}^m [c_i(\phi_i,0)+d_i(0,\phi_i)]$ be an arbitrary element of $\DESH$.
    It follows from \eqref{eq:Comparison}, Lemma~\ref{lem:ChNorm},
    Lemma~\ref{lem:PMRESError}
    and Lemma~\ref{lem:Observations}    that
\begin{align*}
  &\Big\|(\Ch-\Chk)\sum_{i=1}^m [c_i(\phi_i,0)+d_i(0,\phi_i)]\Big\|_\aP\\
  &\hspace{40pt}\leq
  \sum_{i=1}^m \big[|c_i|\|\bpsi_i-\bpsi_{i,k}\|_\aP+|d_i|\|\bxi_i-\bxi_{i,k}\|_\aP\big]\\
  &\hspace{40pt}\leq C_\star q^{\lfloor k/2\rfloor}
  \sum_{i=1}^m \big[|c_i|\|\bpsi_i\|_\aP+|d_i|\|\bxi_i\|_\aP\big]\\
    &\hspace{40pt}\leq C_\star  q^{\lfloor k/2\rfloor} \PFfactor
  \sum_{i=1}^m \big[|c_i|\|\phi_i\|_a+|d_i|\|\phi_i\|_a\big]\\
  &\hspace{40pt}\leq C_\star  q^{\lfloor k/2\rfloor}  \PFfactor \sqrt{\beta}
  \sum_{i=1}^m (|c_i|+|d_i|)|\phi_i|_\HO\\
  &\hspace{40pt}\leq C_\star q^{\lfloor k/2\rfloor}
  \PFfactor  \sqrt{\beta}\, C_\diamondsuit  H^{\tau_d}\sum_{i=1}^m (|c_i|+|d_i|)\\
   &\hspace{40pt}\leq C_\star q^{\lfloor k/2\rfloor}
  \PFfactor  \sqrt{\beta}\, C_\diamondsuit  H^{\tau_d} C_\heartsuit (1/\sqrt{\alpha})H^{-d}
    \Big\|\sum_{i=1}^m \big[c_i(\phi_i,0)+d_i(0,\phi_i)\big]\Big\|_\aP\\
  &\hspace{40pt}= C_\star C_\diamondsuit C_\heartsuit\sqrt{\beta/\alpha}
  \PFfactor H^{-d+\tau_d}q^{\lfloor k/2\rfloor}
     \Big\|\sum_{i=1}^m \big[c_i(\phi_i,0)+d_i(0,\phi_i)\big]\Big\|_\aP,
\end{align*}
 which together with \eqref{eq:C*Def} implies \eqref{eq:OperatorDifference}
 for $C_\spadesuit=2C_\diamondsuit C_\heartsuit$.
\end{proof}
\par
 We have an analog of Corollary~\ref{cor:VHNormEquivalence}.
\begin{corollary}\label{cor:LocalVHNormEquivalence}
  The following relations are valid for any $\bv\in\DESH :$
\begin{align}
 \|\bv-\Chk\bv\|_\aP&\leq \Big(C_\spadesuit
  \kappa(\SKER\AKER)^\frac32 \Big(\frac{\beta}{\alpha}\Big)^\frac12
  \PFfactor^2 H^{-d+\tau_d}q^{\lfloor k/2\rfloor}
   \label{eq:LocalVHNormEquivalence1}\\
    &\hspace{60pt}+  [2+(\CPF/\alpha)]\Big)\|\bv\|_\aP,\notag\\
 \|\bv\|_\aP&\leq C_\dag\sqrt{\beta/\alpha}\|\bv-\Chk\bv\|_\aP.\label{eq:LocalVHNormEquivalence2}
\end{align}
\end{corollary}
\begin{proof}
  The estimate \eqref{eq:LocalVHNormEquivalence1} follows
  from \eqref{eq:VHNormEquivalence1},
  \eqref{eq:OperatorDifference} and the triangle inequality.  The proof of
  \eqref{eq:LocalVHNormEquivalence2}
  is identical to the proof of \eqref{eq:VHNormEquivalence2}.
\end{proof}
\subsection{The Well-posedness of \eqref{eq:LMFEM}}\label{subsec:WP}
 We will use Corollary~\ref{cor:VHNormEquivalence},
   Lemma~\ref{lem:OperatorDifference} and  Corollary~\ref{cor:LocalVHNormEquivalence}
   to derive an analog of Lemma~\ref{lem:Stability} under some assumptions on $k$.
\par\smallskip
\noindent{\bf Assumption 1.} \quad The number of P-MINRES steps
  $k$  is sufficiently large so that, by \eqref{eq:LocalVHNormEquivalence1},
\begin{equation}\label{eq:Assumption1}
  \|\bv-\Chk\bv\|_\aP\leq [3+(\CPF/\alpha)]\|\bv\|_\aP\qquad\forall\,\bv\in\DESH.
\end{equation}
\begin{remark}\label{rem:Assumption1}\rm
  The following  condition on $k$  guarantees \eqref{eq:Assumption1}:
 \begin{equation*}
  C_\spadesuit \Cond^\frac32 \Big(\frac{\beta}{\alpha}\Big)^\frac12
  \PFfactor^2 H^{-d+\tau_d}q^{\lfloor k/2\rfloor}\leq 1,
 \end{equation*}
 or equivalently
\begin{equation}\label{eq:Explicit1}
  \ln C_\spadesuit+\frac32\ln\Cond+\frac12\ln(\beta/\alpha)+2\ln\PFfactor
  +(-d+\tau_d)\ln H+\lfloor k/2\rfloor \ln q \leq 0
\end{equation}
\end{remark}
\par
 It follows from Corollary~\ref{cor:VHNormEquivalence},
 Corollary~\ref{cor:LocalVHNormEquivalence} and \eqref{eq:Assumption1} that
\begin{alignat}{3}
  \|\bv-\Chk\bv\|_\aP&\leq C_\clubsuit\|\bv-\Ch\bv\|_\aP
  &\qquad&\forall\,\bv\in\DESH,  \label{eq:NormComparison2}\\
  \|\bv-\Ch\bv\|_\aP&\leq C_\clubsuit \|\bv-\Chk\bv\|_\aP
   &\qquad&\forall\,\bv\in\DESH,  \label{eq:NormComparison3}
\end{alignat}
where
\begin{equation}\label{eq:CClubDef}
 C_\clubsuit=[3+(\CPF/\alpha)]C_\dag\sqrt{\beta/\alpha}.
\end{equation}
\par\smallskip\noindent{\bf Assumption 2.}\quad The number of P-MINRES steps
 $k$ is sufficiently large so that, by \eqref{eq:VHNormEquivalence2},
  Lemma~\ref{lem:OperatorDifference}
 and \eqref{eq:LocalVHNormEquivalence2},
\begin{equation}\label{eq:Assumption2}
   C_\clubsuit [2+(\CPF/\alpha)]^2\|(\Ch-\Chk)\bv\|_\aP\leq \frac13
    \min(\|\bv-\Ch\bv\|_\aP,\|\bv-\Chk\bv\|_\aP)
\end{equation}
 for all $\bv\in\DESH$.
\begin{remark}\label{rem:Assumption2}\rm
 The following  condition on $k$  guarantees \eqref{eq:Assumption2}:
\begin{equation*}
  3C_\spadesuit C_\dag^2[3+(\CPF/\alpha)]^5\Cond^\frac32\Big(\frac{\beta}{\alpha}\Big)^\frac32
   H^{-d+\tau_d}q^{\lfloor k/2\rfloor}\leq 1,
\end{equation*}
 or equivalently
\begin{align}\label{eq:Explicit2}
  &\ln (3C_\spadesuit C_\dag^2)+5\ln[3+(\CPF/\alpha)]+\frac32\ln\Cond\\
  &\hspace{100pt}+\frac32\ln(\beta/\alpha)+(-d+\tau_d)\ln H+\lfloor k/2\rfloor\ln q\leq 0.\notag
\end{align}
\end{remark}
 \par
  The well-posedness of \eqref{eq:LMFEM} for a sufficiently large
  $k$ is addressed by the following lemma.
\begin{lemma}\label{lem:InfSupLMFEM}
  The inf-sup condition
\begin{equation}\label{eq:InfSupLMFEM}
  \inf_{\strut\bv\in\DESH}\sup_{\substack{\vspace{1pt}\\\bw\in\DESH}}
  \frac{\B{\bv-\Chk\bv}{\bw-\Chk\bw}}{\|\bv-\Chk\bv\|_\aP\|\bw-\Chk\bw\|_\aP}
   \geq \frac{1}{3C_\clubsuit[2+(\CPF/\alpha)]}
\end{equation}
 holds under Assumption $1$ and Assumption $2$.
\end{lemma}
\begin{proof} Let $\bv\in\DESH$ be arbitrary.  We have, by \eqref{eq:Bounded},
 Remark~\ref{rem:MV}, Lemma~\ref{lem:Stability},
 \eqref{eq:NormComparison2}, \eqref{eq:NormComparison3} and \eqref{eq:Assumption2},

\begin{align*}
  \|\bv-\Chk\bv\|_\aP&\leq C_\clubsuit\|\bv-\Ch\bv\|_\aP\\
  &\leq C_\clubsuit [2+(\CPF/\alpha)]\sup_{\bw\in\DESH}
  \frac{\B{\bv-\Ch\bv}{\bw-\Ch\bw}}{\|\bw-\Ch\bw\|_\aP}\\
  &\leq C_\clubsuit [2+(\CPF/\alpha)]\sup_{\bw\in\DESH}
  \frac{\B{(\Chk-\Ch)\bv}{\bw-\Ch\bw}}{\|\bw-\Ch\bw\|_\aP}\\
  &\hspace{20pt}+
    C_\clubsuit [2+(\CPF/\alpha)]\sup_{\bw\in\DESH}
    \frac{\B{\bv-\Chk\bv}{\bw-\Ch\bw}}{\|\bw-\Ch\bw\|_\aP}\\
  &\leq C_\clubsuit [2+(\CPF/\alpha)]^2 \|(\Chk-\Ch)\bv\|_\aP\\
  &\hspace{20pt}+
   C_\clubsuit [2+(\CPF/\alpha)]\sup_{\bw\in\DESH}
   \frac{\B{\bv-\Chk\bv}{(\Chk-\Ch)\bw}}{\|\bw-\Ch\bw\|_\aP}\\
     &\hspace{40pt}+
   C_\clubsuit [2+(\CPF/\alpha)]\sup_{\bw\in\DESH}
   \frac{\B{\bv-\Chk\bv}{\bw-\Chk\bw}}{\|\bw-\Ch\bw\|_\aP}\\
   &\leq \frac13\|\bv-\Chk\bv\|_\aP\\
     &\hspace{20pt}+C_\clubsuit [2+(\CPF/\alpha)]^2\|\bv-\Chk\bv\|_\aP
     \sup_{\bw\in\DESH}\frac{\|(\Ch-\Chk)\bw\|_\aP}{\|\bw-\Ch\bw\|_\aP}\\
     &\hspace{40pt}+ C_\clubsuit [2+(\CPF/\alpha)]\sup_{\bw\in\DESH}
     \frac{\B{\bv-\Chk\bv}{\bw-\Chk\bw}}{\|\bw-\Ch\bw\|_\aP}\\
     &\leq \frac23\|\bv-\Chk\bv\|_\aP
     +C_\clubsuit [2+(\CPF/\alpha)]
     \sup_{\bw\in\DESH}\frac{\B{\bv-\Chk\bv}{\bw-\Chk\bw}}{\|\bw-\Chk\bw\|_\aP},\\
\end{align*}
 which implies \eqref{eq:InfSupLMFEM}.
\end{proof}
\begin{remark}\label{rem:StandardInfSuP}\rm
 In view of \eqref{eq:MSVDef} and \eqref{eq:ChkDef}, we can rewrite the inf-sup condition
 \eqref{eq:InfSupLMFEM} as
\begin{equation}\label{eq:InfSupLMFEM2}
  \inf_{\strut\bv\in\MSV}\sup_{\bw\in\MSV}
  \frac{\B{\bv}{\bw}}{\|\bv\|_\aP\|\bw\|_\aP}
   \geq \frac{1}{3C_\clubsuit[2+(\CPF/\alpha)]}.
\end{equation}
\end{remark}
\subsection{Energy Error}\label{subsec:LMFEMEnergyError}
 Under the inf-sup condition \eqref{eq:InfSupLMFEM2} we have a
 standard quasi-optimal error estimate
 (cf. \cite[Theorem~2.1]{Brezzi:1974:SPP})
 for the solution $\LMSpy$ of \eqref{eq:LMFEM}:
\begin{equation}\label{eq:QuasiOptimal}
  \|\pyh-\LMSpy\|_\aP\leq C_\maltese \inf_{\qz\in\MSV}
  \|\pyh-\qz\|_\aP,
\end{equation}
 where
\begin{equation}\label{eq:CMaltese}
  C_\maltese=1+3C_\clubsuit[2+(\CPF/\alpha)].
\end{equation}
\par
 It then follows from Remark~\ref{rem:AlternativePiHEst},
 \eqref{eq:IMCh}, Remark~\ref{rem:ErrorRepresentation}, \eqref{eq:MSVDef},
  Lemma~\ref{lem:OperatorDifference}
 and \eqref{eq:QuasiOptimal} that
\begin{align}\label{eq:LMFEMEnergyError1}
  &\|\pyh-\LMSpy\|_\aP\notag\\
   &\hspace{30pt}\leq C_\maltese\big\|\pyh-(\bI-\Chk)(\PiH\times\PiH)\pyh\big\|_\aP\notag\\
   &\hspace{30pt}=C_\maltese\big\|\big[\pyh-(\bI-\Ch)(\PiH\times\PiH)\pyh\big]\notag\\
   &\hspace{120pt}+(\Chk-\Ch)(\PiH\times\PiH)\py\big\|_\aP
   \notag\\
   &\hspace{30pt}=C_\maltese\big\|\big[\pyh-(\bI-\Ch)\pyh\big]
   +(\Chk-\Ch)(\PiH\times\PiH)\pyh\big\|_\aP\\
   &\hspace{30pt}\leq C_\maltese\Big(\|\pyh-\MSpy\|_\aP
   +\|(\Chk-\Ch)(\PiH\times\PiH)\pyh\|_\aP\Big)\notag\\
     &\hspace{30pt}\leq C_\maltese\Big(\|\pyh-\MSpy\|_\aP\notag\\
     &\hspace{70pt}+C_\spadesuit C_\dag\Cond^\frac32 ({\beta}/{\alpha})
  \PFfactor^2 H^{-d+\tau_d}q^{\lfloor k/2\rfloor}\|\pyh\|_\aP\Big),\notag
\end{align}
 i.e., up to a term
 that decreases exponentially as $k$ increases, the
 performance of the localized multiscale finite element method defined by \eqref{eq:LMFEM} is
 similar to the performance of the ideal multiscale finite element method defined by
 \eqref{eq:IdealMS}.
\par\smallskip\noindent
{\bf Assumption 3.} \quad The number of P-MINRES steps $k$ is sufficiently large so that
\begin{equation}\label{eq:Assumption3}
  q^{\lfloor k/2\rfloor}\leq \Cond^{-\frac32} (\alpha/\beta) H^{1+d-\tau_d}.
\end{equation}
\begin{theorem}\label{thm:LMFEMEnergyError}
  Under Assumptions~{\rm 1--3}, we have
\begin{equation}\label{eq:LMFEMEnergyError}
  \|\pyh-\LMSpy\|_\aP\leq C_\sharp H\|y_d\|_\LT,
\end{equation}
 where
\begin{equation}\label{eq:CSharp}
  C_\sharp=C_\maltese(C_\dag/\sqrt{\alpha})\Big(1+C_\spadesuit\PFfactor^2\sqrt{\CPF}\Big)
\end{equation}
 is independent of  $h$, $H$ and $\beta/\alpha$.
\end{theorem}
\begin{proof}  The estimate \eqref{eq:LMFEMEnergyError1} is valid under
 Assumption~1 and Assumption~2.
 It then follows from Remark~\ref{rem:WPBound},
 Theorem~\ref{thm:IdealEnergyError} and Assumption~3 that
\begin{align*}
  &\|\pyh-\LMSpy\|_\aP\\
  &\hspace{40pt}\leq C_\maltese\Big((C_\dag/\sqrt{\alpha}H\|y_d\|_\LT
  +C_\spadesuit C_\dag\PFfactor^2 H\|\pyh\|_\aP\Big)\\
  &\hspace{40pt}\leq C_\maltese (C_\dag/\sqrt{\alpha})
  \Big(1+C_\spadesuit\PFfactor^2\sqrt{\CPF}\Big)H\|y_d\|_\LT.
\end{align*}

\end{proof}
\begin{remark}\label{rem:ChoiceOfk}\rm
 Note that \eqref{eq:Assumption3} is equivalent to
\begin{equation}\label{eq:Explicit3}
  \lfloor k/2\rfloor\ln q+\frac32\ln\Cond+\ln(\beta/\alpha)-(1+d-\tau_d)\ln H\leq 0.
\end{equation}
\par
 By examining \eqref{eq:Explicit1}, \eqref{eq:Explicit2} and \eqref{eq:Explicit3},
 we see that the impacts  of the mesh-independent quantities, including the
 condition number $\Cond$ and the contrast $\beta/\alpha$,  are
 mitigated by the natural log  function, and the dominating condition on $k$ is roughly
 (cf. \eqref{eq:Explicit3})
    $$k\geq 2(1+d-\tau_d)\frac{\ln(1/H)}{\ln (1/q)}$$
 From \eqref{eq:qDef} and the relation
     $|\ln(1+x)|\approx |x|$ (for $|x|$ small),
 we can also see that
    $$\ln (1/q)\approx \frac{1}{\Cond}.$$
 Therefore, provided $\Cond$ is moderate, we can choose $k=j\lceil\ln(1/H)\rceil$ for a moderate
 positive integer $j$.
\end{remark}
\subsection{$L_2$ Error}\label{subsec:LMFEMLTwoError}
 We can use a duality argument to obtain an $L_2$ error estimate.
\par
 Observe that Assumption~3 and Lemma~\ref{lem:OperatorDifference} imply
\begin{equation}\label{eq:Assumption4}
  \|(\Ch-\Chk)\bv\|_\aP\leq C_\ddag H\|\bv\|_\aP \qquad\forall\,\bv\in\DESH,
\end{equation}
 where
\begin{equation}\label{eq:Cddag}
 C_\ddag=C_\spadesuit\PFfactor^2.
\end{equation}
\begin{theorem}\label{thm:LMFEMLTwoError}
 Under Assumptions~{\rm 1--3},  we have
\begin{equation}\label{eq:LMFEMLTwoError}
 \|\pyh-\LMSpy\|_\LTP\leq C_\flat\sqrt{\beta/\alpha}\,H^2\|y_d\|_\LT,
\end{equation}
 where
\begin{equation}\label{eq:Cflat}
C_\flat=2C_\dag C_\sharp C_\ddag\PFfactor \sqrt{\CPF/\alpha}
\end{equation}
 is independent of $h$, $H$ and $\beta/\alpha$.
\end{theorem}
\begin{proof}
  Let $\qz\in\DESh$ be defined by
\begin{equation}\label{eq:LocalDuality1}
  \B{\qz}{\rs}=\int_\O (p_h-\LMSp)r\,dx+\int_\O (y_h-\LMSy)s\,dx \qquad \forall\,\rs\in\DESh.
\end{equation}
\par
 It follows from \eqref{eq:Comparison}, \eqref{eq:PF},
 \eqref{eq:InfSup} and \eqref{eq:LocalDuality1} that
\begin{align}\label{eq:WellPosed}
  \|\qz\|_\aP&\leq \|\pyh-\LMSpy\|_\LTP \sup_{\rs\in \ES}\frac{\|\rs\|_\LTP}{\|\rs\|_\aP}\\
   &\leq \sqrt{\CPF/\alpha}\,\|\pyh-\LMSpy\|_\LTP.\notag
\end{align}
\par
 By repeating the arguments in the proof of Theorem~\ref{thm:IdealLTwoError}
 that led to \eqref{eq:Duality2}, we also have
\begin{equation}\label{eq:LocalDuality2}
  \|\Ch \qz\|_\aP\leq (C_\dag/\sqrt{\alpha})H\|\pyh-\LMSpy\|_\LTP.
\end{equation}
\par
 Let $\rs=\pyh-\LMSpy$.  Then the Galerkin relation
\begin{equation}\label{eq:LocalDuality3}
  \B{\bv}{\rs}=0 \qquad \forall\,\bv\in \MSV
\end{equation}
 follows from \eqref{eq:FineDiscreteProblem} and \eqref{eq:LMFEM}.
\par
 Combining \eqref{eq:IMCh}, \eqref{eq:MSVDef},
 \eqref{eq:LocalDuality1} and \eqref{eq:LocalDuality3}, we obtain
\begin{align}\label{eq:LocalDuality4}
  \|\rs\|_\LTP^2&=\B{\qz}{\rs}\notag\\
  &=\B{\qz-(\bI-\Chk)(\PiH\times\PiH)\qz}{\rs}\notag\\
  &=\B{(\Chk-\Ch)(\PiH\times\PiH)\qz}{\rs}\\
  &\hspace{30pt}
  +\B{\qz-(\bI-\Ch)\qz}{\rs}\notag\\
  &=\B{(\Chk-\Ch)(\PiH\times\PiH)\qz}{\rs}+\B{\Ch\qz}{\rs}.\notag
\end{align}
\par
 Using \eqref{eq:Bounded}, Remark~\ref{rem:AlternativePiHEst},
  \eqref{eq:Assumption4}, \eqref{eq:WellPosed},
 we can bound the first term on the right-hand side of \eqref{eq:LocalDuality4} by
\begin{align}\label{eq:LocalDuality5}
  &\B{(\Chk-\Ch)(\PiH\times\PiH)\qz}{\rs}\notag\\
  &\hspace{40pt}\leq
  \PFfactor\big\|(\Chk-\Ch)(\PiH\times\PiH)\qz\big\|_\aP\|\rs\|_\aP\notag\\
  &\hspace{40pt}\leq \PFfactor C_\ddag H\|(\PiH\times\PiH)\qz\big\|_\aP\|\rs\|_\aP\\
  &\hspace{40pt}\leq \PFfactor C_\ddag  C_\dag\sqrt{\beta/\alpha}H\|\qz\|_\aP\|\rs\|_\aP\notag\\
  &\hspace{40pt}\leq \PFfactor  C_\ddag  C_\dag\sqrt{\beta/\alpha}
  \sqrt{\CPF/\alpha}H\|\rs\|_\LTP\|\rs\|_\aP.\notag
\end{align}
\par
 For the second term on the right-hand side of \eqref{eq:LocalDuality4}, we have the bound
\begin{align}\label{eq:LocalDuality6}
  \B{\Ch\qz}{\rs}
  &\leq
  \PFfactor\|\Ch\qz\|_\aP\|\rs\|_\aP\\
  &\leq \PFfactor (C_\dag/\sqrt{\alpha})H\|\rs\|_\LTP\|\rs\|_\aP\notag
\end{align}
 by \eqref{eq:Bounded} and \eqref{eq:LocalDuality2}.
\par
 Putting \eqref{eq:LMFEMEnergyError} and \eqref{eq:LocalDuality4}--\eqref{eq:LocalDuality6} together,
 we arrive at the estimate
\begin{align*}
 &\|\pyh-\LMSpy\|_\LTP\\
 &\hspace{30pt}\leq \PFfactor (C_\dag/\sqrt{\alpha})\big[C_\ddag\sqrt{\CPF(\beta/\alpha)}+1\big]H
 \|\pyh-\LMSpy\|_\aP\\
 &\hspace{30pt}\leq \PFfactor (C_\dag/\sqrt{\alpha})\big[C_\ddag\sqrt{\CPF(\beta/\alpha)}+1\big]C_\sharp
   H^2\|y_d\|_\LT,
\end{align*}
 which implies \eqref{eq:LMFEMLTwoError} with $C_\flat$ given by \eqref{eq:Cflat}.
\end{proof}
\section{Numerical Results}\label{sec:Numerics}
 The numerical results for Example~\ref{example:Oscillatory} and
 Example~\ref{example:Heterogeneous} are presented in Section~\ref{subsec:HOP} and
 Section~\ref{subsec:HHP}.  We also describe briefly some computational aspects in
 Section~\ref{subsec:HPC}.
\subsection{Highly Oscillatory Problem}\label{subsec:HOP}
  We solve the optimal control problem \eqref{eq:OCP} on the unit square $(0,1)\times(0,1)$,
  where
  $\gamma=1$, $y_d=-1$, and $\cA$ is the matrix in Example~\ref{example:Oscillatory}
  with $\epsilon=0.08$, $0.04$ and $0.025$.  We use the localized multiscale finite element method from
   Section~\ref{sec:LMFEM}.
 \par
 For this problem $\alpha\approx 1$ and $\beta\approx 20$, the magnitude of
 $\Cond$ is moderate by \eqref{eq:qBound} and, according to
 Remark~\ref{rem:ChoiceOfk},  we can choose
 the number of P-MINRES steps $k$ to be $j\lceil\log(1/H)\rceil$ for a moderate positive integer $j$.
   The choices for $h$, $H$ and $j$ for different values of $\epsilon$
   are described in Table~\ref{table:OscillatoryData}.
\par
\begin{table}[h]
\begin{tabular}{|c|c|c|c|c||c|c|c|c|}
  \hline
  $\epsilon$ &\multicolumn{4}{c||}{$0.08$\quad$0.04$}&\multicolumn{4}{c|}{$0.025$}\\
  \hline
  $1/h$ &\multicolumn{4}{c||}{$256$}& \multicolumn{4}{c|}{$320$}\\
  \hline
  $1/H$ & $8$ & $16$ & $32$ & $64$ & $10$ & $20$ & $40$ & $80$ \\
  \hline
  $j$ & $2$ & $2$ & $3$ & $3$ & $2$ & $2$ & $3$ & $3$ \\
  \hline
\end{tabular}
\par\medskip
\caption{Choices for $h$, $H$, and $j$ for different values of $\epsilon$.}
\label{table:OscillatoryData}
\end{table}
\par
\begin{figure}[htb!]
  \centering
  \includegraphics[width=0.6\linewidth]{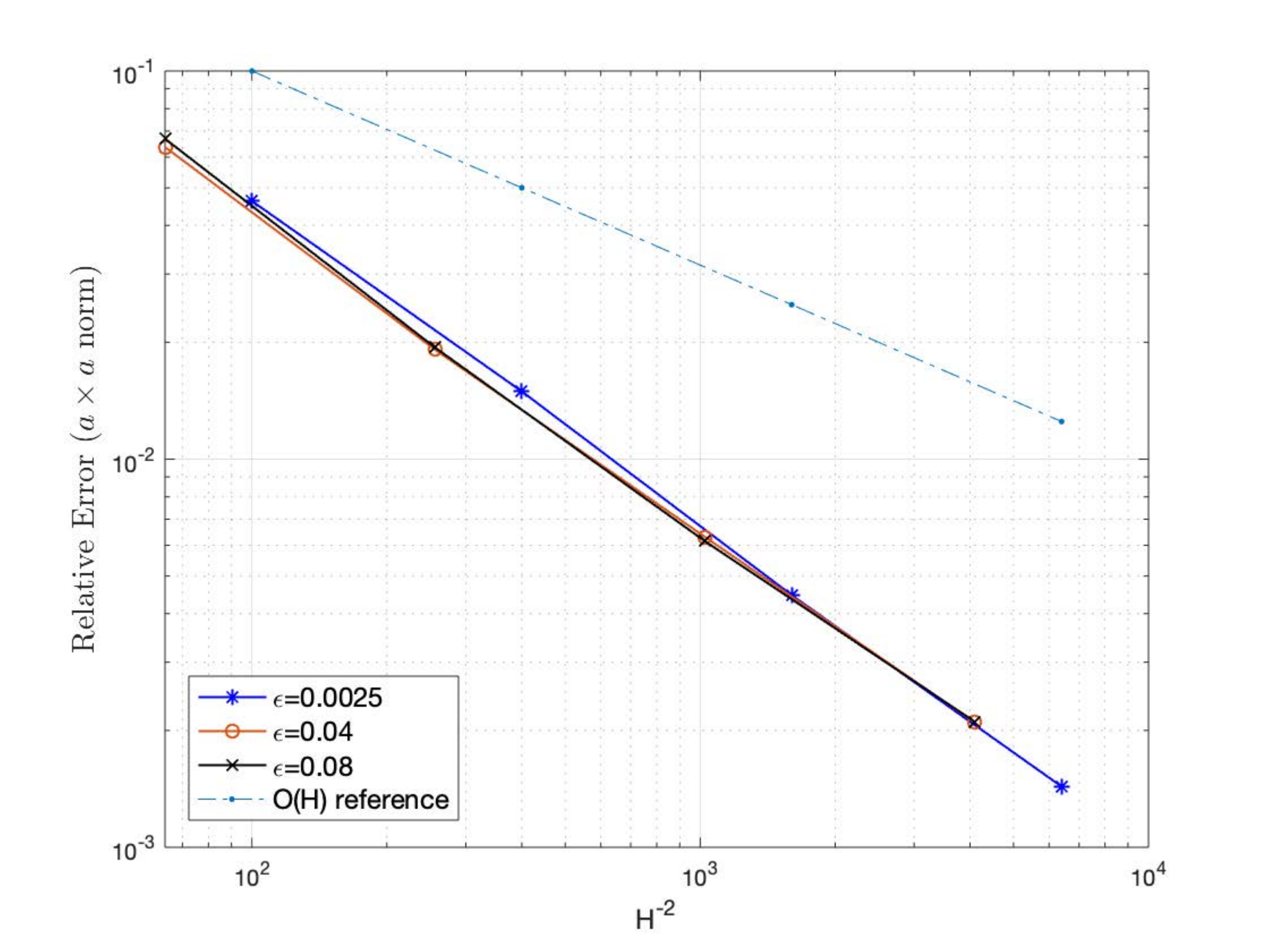}  
  \caption{Relative $\|\cdot\|_\aP$  errors of the localized multiscale approximations
  of the highly oscillatory problem for different values of $\epsilon$.}
  \label{fig:OscillatoryH1}
\end{figure}
\begin{figure}[htb!]
  \centering
  \includegraphics[width=0.6\linewidth]{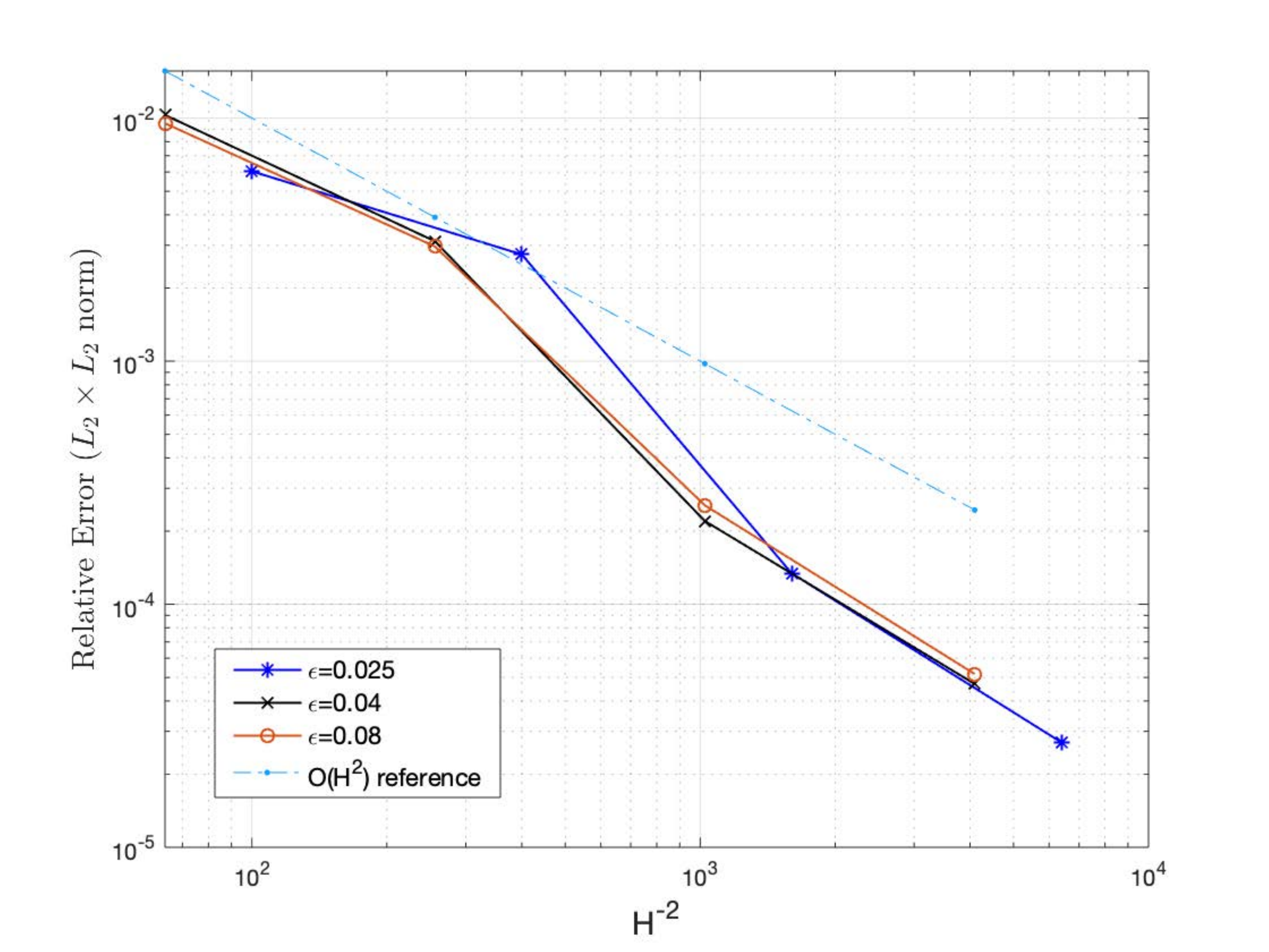}
  \caption{Relative $\|\cdot\|_\LTP$  errors of the localized multiscale approximations
  of the highly oscillatory problem for different values of $\epsilon$ }
  \label{fig:OscillatoryL2}
\end{figure}

 The relative errors in the $\|\cdot\|_\aP$ norm and the $\|\cdot\|_\LTP$ norm are
 presented in Figure~\ref{fig:OscillatoryH1}
 and Figure~\ref{fig:OscillatoryL2}.  We observe $O(H)$ convergence in the $\|\cdot\|_\aP$ norm
 and $O(H^2)$ convergence in the $\|\cdot\|_\LTP$ norm,
 which agree with
 Theorem~\ref{thm:LMFEMEnergyError} and Theorem~\ref{thm:LMFEMLTwoError}.
\subsection{Highly Heterogeneous Problem}\label{subsec:HHP}
  We solve the optimal control problem \eqref{eq:OCP} on the unit square, where $\gamma=1$, $y_d=1$, and
  $\cA$ is the matrix from Example~\ref{example:Heterogeneous}.  We use the ideal multiscale finite element method
  from Section~\ref{sec:Ideal} and the localized multiscale finite element method from
  Section~\ref{sec:LMFEM}.  The reference mesh size is $h=1/320$.
\par
  For this problem we have $\alpha=1$ and $\beta=1350$. The value of the condition number
  $\Cond$  is found computationally to be less than $10$, which is much better than the pessimistic
  bound in \eqref{eq:ConditionNumber}.  Therefore we can, according to Remark~\ref{rem:ChoiceOfk},
  choose the number of P-MINRES steps $k$ to be
  $j\lceil\log(1/H)\rceil$ for a moderate positive integer $j$.
  Here we take $j=2$ for $H=1/10$, $j=3$ for $H=1/20$, and $j=4$ for $H=1/40$.
\par
\begin{figure}[htb!]
  \centering
  \includegraphics[width=0.6\linewidth]{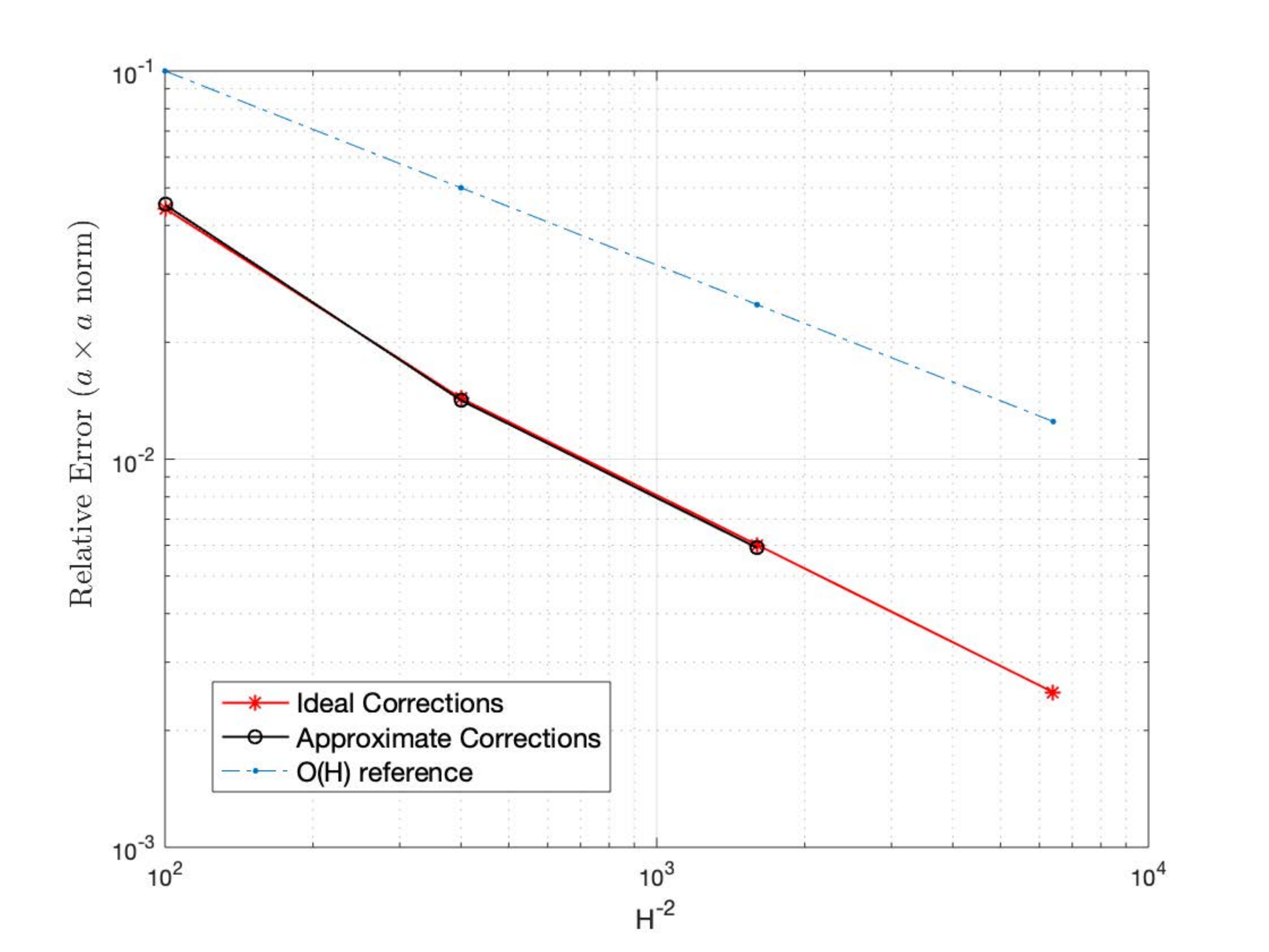}
  \par\medskip
 \caption{Relative $\|\cdot\|_\aP$  errors of the multiscale approximations
  of the highly heterogeneous problem,
   where $h=1/320$,
   $k=j\lceil\log(1/H)\rceil$, with $j=2$ for $H=1/10$, $j=3$ for $H=1/20$, and $j=4$ for $H=1/40$.}
  \label{fig:HeterogeneousH1}
\end{figure}
\begin{figure}[htb!]
  \centering
  \includegraphics[width=0.6\linewidth]{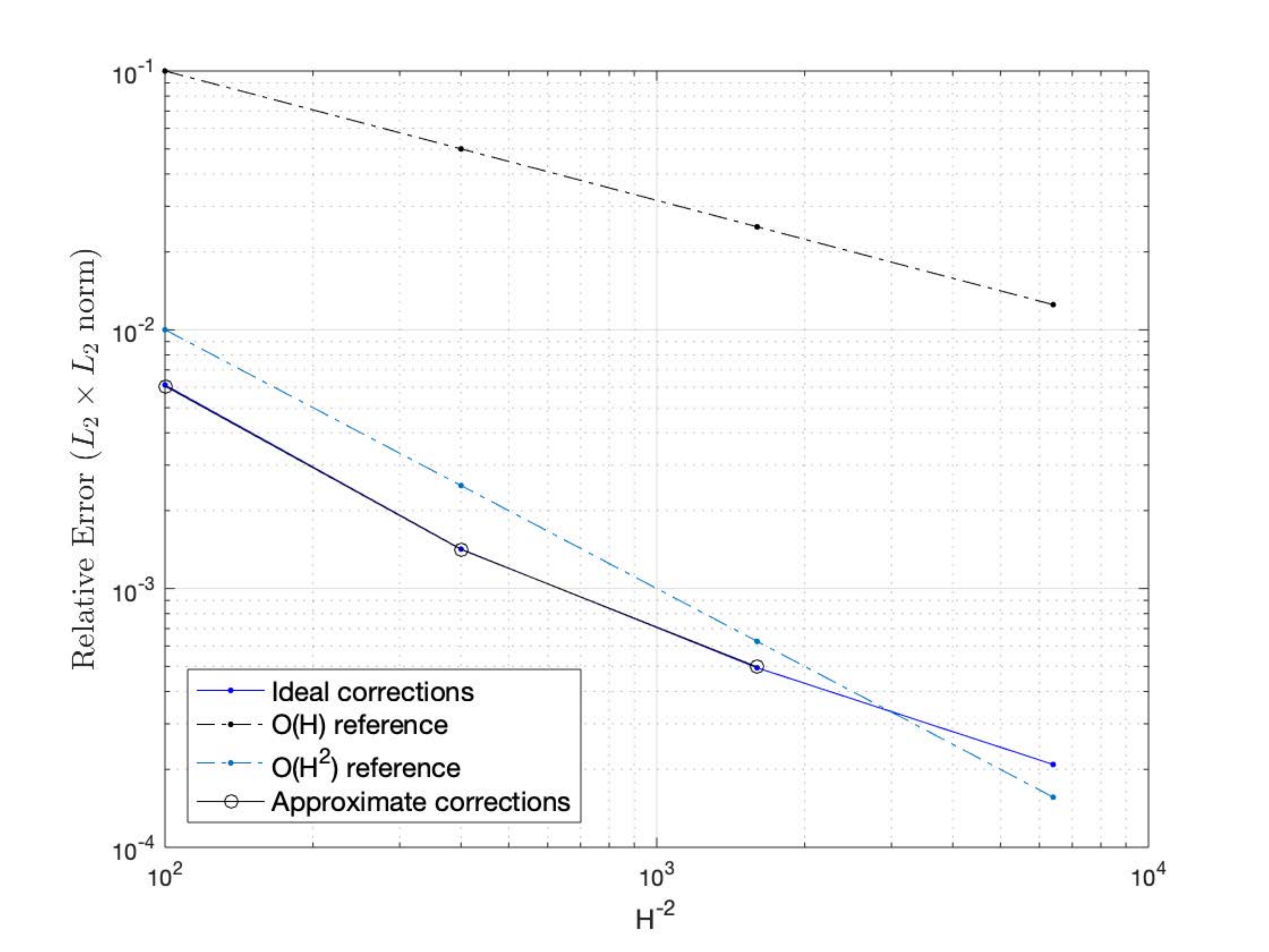}
  \caption{Relative $\|\cdot\|_\LTP$  errors of multiscale approximations
  of the highly heterogeneous problem,
   where $h=1/320$,
   $k=j\lceil\log(1/H)\rceil$, with $j=2$ for $H=1/10$, $j=3$ for $H=1/20$, and $J=4$ for
   $H=1/40$.}
  \label{fig:HeterogeneousL2}
\end{figure}

 The relative errors in the $\|\cdot\|_\aP$ norm and the $\|\cdot\|_\LTP$ norm are displayed in
 Figure~\ref{fig:HeterogeneousH1} and Figure~\ref{fig:HeterogeneousL2}.  We observe that the errors for the ideal
 multiscale finite element method and the localized multiscale finite element method are indistinguishable.
 The order of convergence in the $\|\cdot\|_\aP$ norm is 1, which agrees with Theorem~\ref{thm:IdealEnergyError}
 and Theorem~\ref{thm:LMFEMEnergyError}.  The convergence history for the $\|\cdot\|_\LTP$ norm is similar to
 the early stage of the history for the highly oscillatory but well-conditioned problem in Figure~\ref{fig:OscillatoryL2}.  Therefore it is reasonable to expect that the order of convergence
 in the $\|\cdot\|_\LTP$ norm for the ill-conditioned highly heterogenous problem will also approach
 $2$ at a later stage.
%
\subsection{Some Computational Aspects}\label{subsec:HPC}
 We will focus on  the highly heterogeneous problem in Section~\ref{subsec:HHP}, where the reference
 solution is obtained by a standard finite element method with $h=1/320$.  Below are some observations on
 the case where the coarse mesh size is $H=1/20$ and the solution obtained by the localized multiscale method in
 Section~\ref{sec:LMFEM} is quite reasonable (cf. Figure~\ref{fig:HeterogenousFEMs}).
\par
 The parallel computing was carried out on a cluster with 440 compute nodes  running the Red Hat Enterprise Linux 6 operating system and has a 146 TFlops peak performance.
 Each compute node is equipped with two 8-core Sandy Bridge Xeon 64-bit processors operating at a core frequency of
 2.6 GHz, 32GB 1666MHz RAM, 500GB HD, 40 Gigabit/sec Infiniband network interface and a 1 Gigabit Ethernet network interface.
\par
 We computed the reference solution  with the PETSc library
 using  128 processors and an ILU(0) preconditioner in PGMRES.
   The solution time is $0.22$ seconds.
 For comparison, we solved the smaller system \eqref{eq:LMFEM} by Gaussian elimination in MATLAB on a MacBook Pro
  ($2.8$ Ghz Quad-Core Intel Core i$7$ processor and a $16$GB $2133$ Mhz LPDDR RAM).
 The solution time is $0.02$ seconds.
\par
 The total (set-up and solution) time for computing the reference solution with PETSc using 128
 processors is $0.47$ seconds.  For comparison, the total time for solving  \eqref{eq:LMFEM}
 with 128 different right-hand sides simultaneously
 using PETSc and Gaussian elimination is $1.36$ seconds.
\par
  Using 1024 processors, it took 2104 seconds in the offline stage to construct the basis functions
 of the ideal multiscale
 finite element space $\MV$ in Section~\ref{subsec:MV}, and 535 seconds to construct the basis functions
  of the localized  multiscale
 finite element space $\MSV$ in Section~\ref{subsec:MSV}.
%
\section{Concluding Remarks}\label{sec:Conclusion}
 In this paper we have developed  multiscale finite element methods for a linear-quadratic
  elliptic optimal control problem with rough coefficients, where scale separation and periodic structures
  are not assumed.   These methods can be viewed as reduced order methods.
 \par
  In particular, we have constructed
  a generalized finite element method with localized basis functions whose performance is
  similar to standard finite element methods for smooth problems.
  Both the construction of the generalized finite element space
  and the analysis of the resulting Galerkin method are
  based on basic finite element technology and two (by now) classical numerical linear algebra
  ingredients, namely  the additive Schwarz preconditioner and the preconditioned minimum
  residual algorithm.  Our work further illustrates the idea put forth in \cite{KPY:2018:LOD}
  that multiscale problems can be solved by combining finite element methods, domain decomposition algorithms
   and   iterative Krylov subspace solvers.
\par
  The techniques developed in this paper and \cite{BGS:2021:LOD} can also be extended to elliptic variational
  inequalities with rough coefficients, such as the obstacle problem and the optimal control problem with
  control constraints.
\section*{Acknowledgements}
Portions of this research were conducted with high performance computing
resources provided by Louisiana State University (http://www.hpc.lsu.edu).
\section*{Data Availability}
The datasets generated during and/or analysed during the current study are available from the corresponding author on reasonable request.

\end{document}